\documentclass[twocolumn]{autart}
\usepackage[namelimits,reqno]{amsmath}
\usepackage{mathptmx} 
\usepackage{amsmath} 
\usepackage{amssymb}
\usepackage{amsfonts}
\usepackage{graphicx}
\usepackage{subfigure}
\usepackage{epstopdf}
\usepackage{bm}
\usepackage{float}
\usepackage{latexsym}
\usepackage{epstopdf}
\usepackage[english]{babel}
\usepackage{epsfig}
\usepackage{fancyhdr}
\usepackage{eepic}
\usepackage{pstricks}
\usepackage{textpos}
\usepackage{niceframe}
\usepackage{fancybox}
\usepackage{hyperref}
\usepackage{tikz} \usetikzlibrary{calc} \tikzset{>=latex}

\newcommand{\beqn}{\begin{align}}
\newcommand{\eeqn}{\end{align}}

\newcommand{\beqnn}{\begin{equation*}}
\newcommand{\eeqnn}{\end{equation*}}

\newcommand{\bsubeq}{\begin{subequations}}
\newcommand{\esubeq}{\end{subequations}}

\newcommand{\bma}{\left(\begin{array}}
\newcommand{\ema}{\end{array}\right)}

\newcommand{\biz}{\begin{itemize}}
\newcommand{\eiz}{\end{itemize}}

\newcommand{\benu}{\begin{enumerate}}
\newcommand{\eenu}{\end{enumerate}}

\newcommand{\bce}{\begin{center}}
\newcommand{\ece}{\end{center}}

\newcommand{\bem}{\begin{em}}
\newcommand{\eem}{\end{em}}

\newcommand{\bpm}{\begin{pmatrix}}
\newcommand{\epm}{\end{pmatrix}}

\newcommand{\p}{\partial}

\begin{document}

\graphicspath{{picfroude/}}

\newtheorem{propo}{Proposition}[section]
\newtheorem{exple}{Exemple}[section]
\newtheorem{lemm}{Lemma}
\newtheorem{remm}{Remark}

\begin{frontmatter}

\title{Backstepping
Stabilization of\\ the Linearized \textit{Saint-Venant-Exner} Model\thanksref{footnoteinfo}} 

\thanks[footnoteinfo]{Corresponding author: A. Diagne.}

\author[AD]{Ababacar Diagne}\ead{ababacar.diagne@it.uu.se},    
\author[MD]{Mamadou Diagne}\ead{mdiagne@ucsd.edu},               
\author[MD]{Shuxia Tang}\ead{sht015@ucsd.edu},  
\author[MD]{Miroslav Krstic}\ead{krstic@ucsd.edu}

\address[AD]{Division of Scientific Computing, Department of Information Technology,  Uppsala University,
Box 337, 75105 Uppsala, Sweden}  
\address[MD]{Department of Mechanical $\&$ Aerospace Engineering,
University of California, San Diego,
La Jolla, CA 92093-0411}             


\begin{abstract}
Using the backstepping design, we achieve exponential stabilization of the coupled Saint-Venant-Exner (SVE) PDE model of water dynamics
in a sediment-filled canal with arbitrary values of canal bottom slope, friction, porosity, and water-sediment  interaction under subcritical or supercritical flow regime.
The studied  SVE model consists of two rightward convecting transport Partial Differential Equations (PDEs) and one leftward convecting transport
PDE. A single boundary input control (with actuation located only at downstream) strategy is  adopted. 
A full state feedback controller is firstly designed, which guarantees the exponential stability  of the closed-loop control system.
Then, an  output feedback controller is designed based on  the reconstruction of the distributed  state with  a
backstepping observer. It also guarantees the exponential stability  of the closed-loop control system. 
The flow regime depends on the dimensionless Froude number $Fr$, and both our controllers can deal 
with the subcritical ($Fr<1$) and supercritical ($Fr>1$) flow regime. They achieve the exponential stability results   
without  any restrictive conditions in contrast to existing results.
\end{abstract}

\end{frontmatter}

\section{Introduction}
Balance laws are the key point for modeling complex physical systems that involve  fluid mechanics, reactions, heat and  mass transfer phenomena.  In  fluid mechanics, fundamental  balance equations   expressing the conservation   of certain quantities, such as the energy, the mass or the momentum in physical processes, lead to spatio-temporal differential equations that  express transport or  diffusion phenomena. Such equations are the starting point for the design of various controllers that ensure   the stability and operability  of many engineering applications.  Among those applications, we are interested in the stabilization of the hyperbolic  SVE PDEs describing   the flow and the  bed evolutions in an open channel \cite{Diagne:2011,Diagne:2012}.  The SVE model has attracted considerable attention over the
past decades. Several  theoretical and experimental studies have been proposed in the literature, considering  the flow and sediment characteristics of the water motion. These studies also addressed the influence of the  particle size, shape and density.  However, the control of such systems, modeled by nonlinear hyperbolic PDEs, is left out in most of the studies. To the best of the authors' knowledge,  there are only a few results on the stabilization of SVE  model in the existing literatures.

 Several strategies  have been developed  to control  the flow dynamics  in classical  irrigation canals without the sediment layer during the last decades. We refer the reader to     \cite{malat} in which the  classification of control problems and related methodologies   is fairly addressed. Basically, the main  purpose is the regulation of  the water level at a desired height  by adjusting the opening of the gates at the ends  of the channel as boundary actuators. For instance, the synthesis of LQ controller can be found in \cite{bhv,malapilot,weyer1}, whereas  \cite{litfr,pognant}   has studied an $H_{\infty}$ control approach. Through semigroup approach, \cite{xu}  proposed an integral output feedback controller using a linear PDE model around a steady state.  Lyapunov analysis is investigated  \cite{CAB2,Prieur2004,DosSantos08}  and  multi-models approach with a stability analysis based on  Linear Matrix Inequality ($LMI$) is presented in  \cite{Mdiagne10,Mdiagne12}.
 
 Very recently,  \cite{tang2014boundary} proposed a  singular perturbation approach for the synthesis of  boundary control  for hyperbolic systems. The effectiveness of the controller is illustrated using the linearized SVE model. In     \cite{Diagne:2011},  explicit boundary dissipative conditions are given for the exponential stability in $L^2$-norm of one dimensional linear hyperbolic systems of balance laws. The proposed Lyapunov approach is applied to the linearized SVE equations  with successful results. However,  on-line measurements of the water levels at both ends of the spatial domain, namely,  at upstream $ (x = 0)$   and downstream $(x = L)$  are assumed to be available. Later on,    a priori estimation technique and the Faedo-Galerkin method
are  proposed in \cite{Diagne:2012}  for the  design of  a linear feedback control law  that requires only downstream measurements.
 Such an approach have been vividly presented in \cite{ben2013,sam2013} as well.

The stabilization problems for hyperbolic systems have been widely studied in the literature.
The first approach  relies  on careful analysis of the classical solutions along
 the characteristics.  We refer the readers to \cite{LiJMG} in the
case of second-order system of conservation laws and to \cite{LiSt1} for  more general situations as  $n$th-order systems.
Another approach based on the Lyapunov techniques is  introduced by  \cite{CAB1} and
  improved in \cite{CAB2} where  a strict Lyapunov function in terms of Riemann invariants
is  constructed and  its time derivative can be made  negative definite by choosing properly the boundary conditions.
The aforementioned Lyapunov function is very useful to analyze nonlinear hyperbolic systems of conservation laws because of its robustness. We refer to
\cite{CoronBook,Coron08,CAB2,CWang13,Diagne:2012,XuCZ09,XuCZ02}, in which several  applications are analyzed based on this tool.

Recently, the backstepping method was introduced for the feedback stabilization of various  classes  of PDEs \cite{Dimeglio2013,krstic2008boundary}. The key idea of this  approach  is the construction of suitable Volterra integral transformations  that map  the original system into a so-called ``target system'', which is exponentially stable. The kernel functions of the transformations are required to satisfy some PDEs, and the solutions can be then used as gains of the controllers. The invertibility of the  transformations  ensures the exponential stability of the closed-loop control systems. One can refer to  \cite{coron2013local,krstic2008boundary,andrew2004,Bernard2014} for further applications  of this technique to other classes of systems   including nonlinear PDEs.

In the present work, we achieve an exponential stabilization of the coupled Saint-Venant-Exner (SVE)  model  of water dynamics in a sediment-filled canal with arbitrary values of canal bottom slope, friction, porosity, and water-sediment interaction \cite{Diagne:2011}.  The studied  SVE model consists of two rightward  and one leftward convecting transport PDEs \cite{Dimeglio2013}. 
First, a  single boundary input control strategy (with actuation located only at downstream) is  adopted and     a full state feedback backstepping controller    is designed  to ensure the exponential stabilization of the system. Next,  we employ a sole sensor at the upstream  $(x=0)$ to derive an output feedback controller based on  the reconstruction of the distributed  state with an exponentially convergent Luenberger observer. The  flow characteristics  depend  on the dimensionless Froude number, namely,  $Fr$, and  the proposed  controllers operate under  the subcritical ($Fr<1$) and the supercritical ($Fr>1$) flow regime. Moreover,  the exponential stability results   are obtained without imposing any \textcolor{black}{restrictive  conditions on the controller gains}, which is in contrast to  \cite{Diagne:2011}.  
In this context, this paper is an  extension of    \cite{Diagne:2011} in which the design of the controller that guarantees the exponential stability property,  requires a measurement at the two boundary of the channel while in the present work we just use the downstream gate for actuation without need of a sensor there.
Better, with the backstepping method,  we achieve the exponential stability around the origin without imposing any conditions
on the matrix from the source term of the system under study.
Conversely in \cite{Diagne:2011},  that matrix is required to satisfy a restrictive condition to be marginally diagonally stable.

This paper is  organized as follows. In  the next Section,  the nonlinear SVE model  is formulated based on its physical description,  and a linearized version around a steady state is presented.   Section \ref{backstepping} is dedicated to the backstepping transformation between the linearized model and a suitable exponentially stable target system. Then, with the solutions to the gain  kernel PDEs of  the Volterra transformation, a full state controller is computed. An exponentially convergent backstepping observer is designed in \ref{observer}. Based on the observer, which reconstructs the full state from the output measurement, an output feedback controller is constructed in section \ref{output} and an exponential stability  result is also established. Numerical simulations are  provided for both the subcritical and supercritical flow regime in Section \ref{simulations}, with a detailed discussion on the numerical computation of the controller gain. Finally in Section \ref{conclusion}, a  conclusion is presented  and some perspectives  are discussed.

\section{The Saint-Venant-Exner model}\label{model}

\begin{figure}[H]
	\centering
	\includegraphics[width=2.3in]{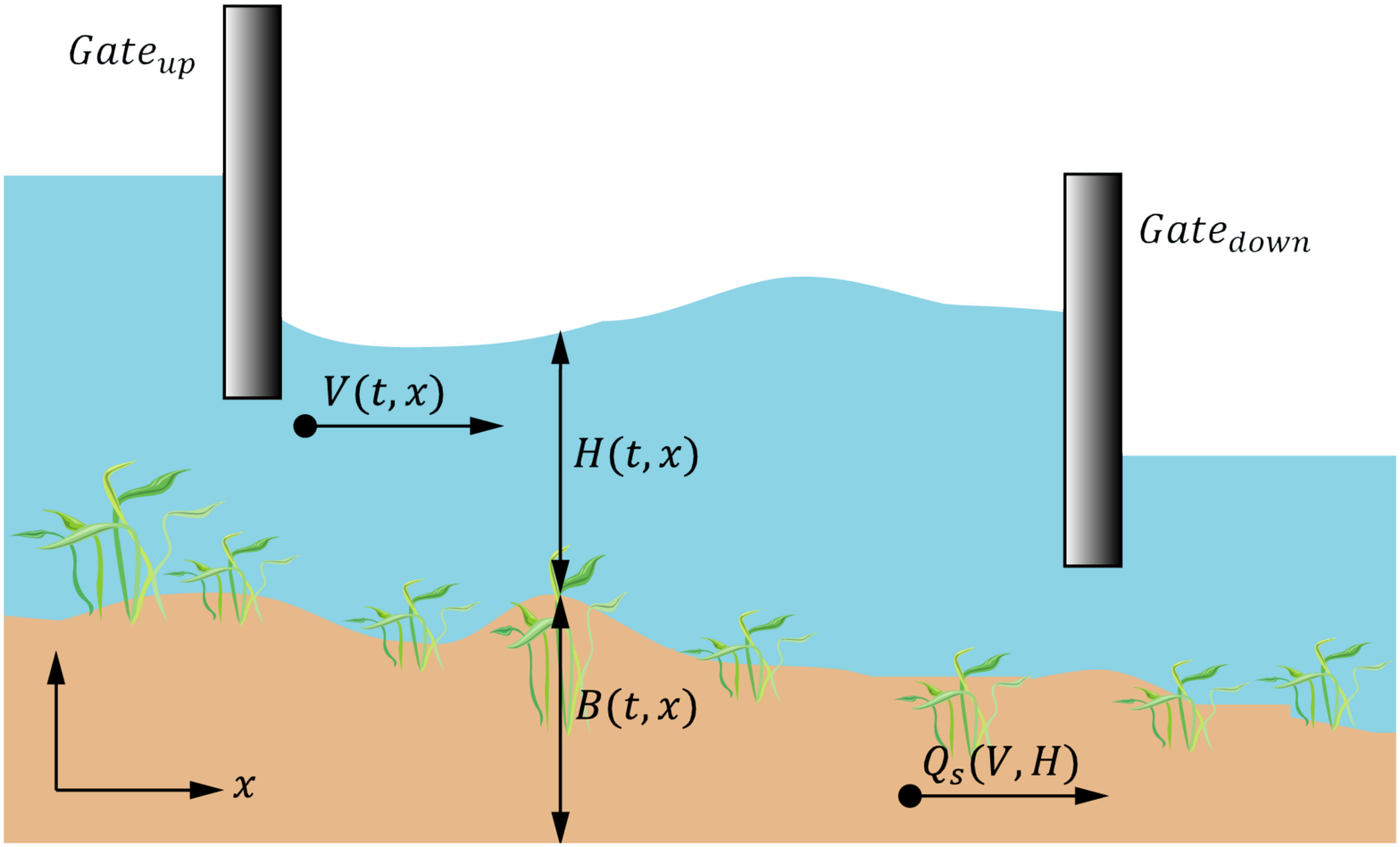}
	 \caption{ A sketch of the channel.}
	  \label{fig:schemacanal}
\end{figure}

We consider a pool of a prismatic sloping open channel with a rectangular cross-section, a unit width and a moving bathymetry
(because of sediment transportation). The state variables of the model are: the water depth $H(t,x)$, the water velocity $V(t,x)$
and the bathymetry $B(t,x)$ which is the depth of the sediment layer above the channel bottom as depicted on Figure  \ref{fig:schemacanal}.
The dynamics of the system are described by the coupling of Saint-Venant and Exner equations (see e.g. \cite{HudSwe03}):
\bsubeq\label{eq1}
\begin{align}
&\dfrac{\partial H}{\partial t} +  V\frac{\partial H}{\partial x} + H\frac{\partial V}{\partial x} \;=\;0 \\
&\dfrac{\partial V}{\partial t} + V\frac{\partial V}{\partial x} + g\frac{\partial H}{\partial x}
+ g\frac{\partial B}{\partial x}\;=\;gS_{b}-C_{f}\frac{V^2}{H} \\
&\dfrac{\partial B}{\partial t} + a V^{2} \frac{\partial V}{\partial x}\;=\; 0.
\end{align}
\esubeq
In these equations, $g$ is the gravity constant, $S_b$ is the bottom slope of the channel, $C_f$ is a friction coefficient and $a$ is
a parameter that encompasses the porosity and viscosity effects on the sediment dynamics. The coefficient $a$ expresses as (cf \cite{HudSwe03})
$$a=\frac{3A_g}{1-p_g}$$
where $p_g$ is the porosity parameter and $A_g$ is the coefficient to control the interaction between the bed and the water flow.
\subsection{Steady-state and Linearization}

A {\it steady-state} is a constant state ($H^*$, $V^*$, $B^*)^T$ which satisfies the relation
\beqnn
gS_b H^* = C_f V^{*2}.
\eeqnn
In order to linearize the model, we define the deviation of the state $\big(H(t,x), V(t,x)$,
$B(t,x)\big)^T$ with respect to the steady-state:
\begin{align*}
&\left(
  \begin{array}{c}
    h(x,t)\\  u(x,t) \\ b(x,t)
  \end{array}
\right)
=
\left(
  \begin{array}{c}
    H(x,t) - H^*\\
     V(x,t) - V^*\\
      B(x,t) - B^*
  \end{array}
\right).
\end{align*}
Then the linearized system of the SVE model (\ref{eq1}) around a steady-state is
\bsubeq \label{eql_1}
\begin{align}
&\dfrac{\partial h}{\partial t}  + V^*\frac{\partial h}{\partial x} + H^*\frac{\partial {u}}{\partial x} =0\\
&\dfrac{\partial {u}}{\partial t} +  V^*\frac{\partial {u}}{\partial x} + g\frac{\partial h}{\partial x}
+ g\frac{\partial b}{\partial x}  =C_{f}\frac{V^{*2}}{H^{*2}}h-2C_{f}\frac{V^*}{H^*}u\\
&\dfrac{\partial b}{\partial t} +  a V^{*2} \frac{\partial {u}}{\partial x}=0.
\end{align}
\esubeq
\subsection{Characteristic (Riemann)  coordinates}
In the matrix form, the linearized model (\ref{eql_1}) can be written as
\begin{align}
 \displaystyle{\frac{\partial W}{\partial t}+ \bm{A}(W^*)\frac{\partial W}{\partial x}=\bm{B}(W^*)W},
\label{ncmfl}
\end{align}
where
\beqnn
W=
\left(
\begin{array}{c}
 {h} \\
 {u} \\
 {b} \\
\end{array}
\right),~
\bm{A}(W^*)=\left(
\begin{array}{ccc}
 {V^*} &  {H^*}  &  {0}\\
 {g} &  {V^*}  &  {g}\\
 {0} &  \displaystyle aV^{*2}  &  {0}
\end{array}
\right),
\eeqnn
\beqnn
\bm{B}(W^*)=\left(
\begin{array}{ccc}
 {0} &  {0}  &  {0}\\
 \displaystyle{C_{f}\frac{V^{*2}}{H^{*2}}} &  \displaystyle{-2C_{f}\frac{V^*}{H^*}}  &  {0}\\
 {0} &  {0}  &  {0}
\end{array}
\right).
\eeqnn
The  dimensionless Froude number is defined as
$$Fr=\frac{V^*}{\sqrt{gH*}}.$$
Exact, but rather complicated expressions of the eigenvalues of $\bm{A}(W^*)$ can be obtained by
using the \textit{Cardano-Vieta} method, see \cite{HudSwe03}.
Once the eigenvalues $\lambda_i$ of the matrix $\bm{A}(W^*)$ are obtained, the corresponding left
eigenvectors can be computed as
\begin{align}
L_k &= {\frac{1}{(\lambda_k-\lambda_i)(\lambda_k-\lambda_j)}} \left(
\begin{array}{c}
 {(V^*-\lambda_i)(V^*-\lambda_j)+gH^*} \\
 {H^*\lambda_k} \\
 {gH^*} \\
\end{array}
\right),  \nonumber\\
&{\rm{for}}~ k \neq i \neq j \;\in\;\left\{1,2,3\right\}.
\label{lefteig}
\end{align}
We  multiply  (\ref{ncmfl}) by $L^{T}_{k}$ in order to rewrite the model in terms of the characteristic
coordinates $\psi_k$ ($k=1,2,3$). Then we obtain
\begin{align}
\frac{\partial \Phi_{k}}{\partial t}+
\lambda_{k}\frac{\partial \Phi_{k}}{\partial x}=L^{T}_{k} \bm{B} W, {\rm{for}}~ k=1,2,3,
\label{rimequat1}
\end{align}
where
\begin{align}\label{eq:psi-k}
\Phi_k =& \dfrac{1}{(\lambda_k-\lambda_i)(\lambda_k-\lambda_j)}\left[ \big((V^*-\lambda_i)(V^*-\lambda_j)\right. \nonumber\\
& \left.
+ gH^*\big) h   + H^*\lambda_k u + gH^*b \right].
\end{align}
 For the sake of simplicity, we introduce the following notation $r_{k}$:
$$r_{k}=C_{f}\frac{V^*}{H^*}\frac{\lambda_k}{(\lambda_k-\lambda_i)(\lambda_k-\lambda_j)}.$$
Some computations yields the following writing for equation (\ref{rimequat1}):
\begin{align}
\dfrac{\partial \xi_{k}}{\partial t}+ \lambda_{k}\frac{\partial \xi_{k}}{\partial x}
+\sum_{s=1}^{3}(2\lambda_{s}-&3V^*)r_{s}\xi_{s}=0,\nonumber\\
& {\text{for}}~ k=1,2,3,
\label{rimequat2}
\end{align}
where the characteristic coordinates are now defined as
\begin{align}
\displaystyle{\xi_{k}=\frac{1}{r_{k}}\Phi_{k}}.
\label{eq:psi-to-xi}
\end{align}
From (\ref{rimequat2}), the linearized model (\ref{rimequat1}) in  characteristic form can be written as
\begin{align}
\frac{\p \bm{\xi}}{\p t}+ \bm{\Lambda}\frac{\p \bm{\xi}}{\p x} - \bm{M}\bm{\xi}=0,
\label{rimgen}
\end{align}
where
\beqnn
\bm{\xi}=(\xi_{1},\xi_{2},\xi_{3})^{T},\quad
\bm{\Lambda}=diag(\lambda_{1},\lambda_{2},\lambda_{3}),
\eeqnn

\beqnn
\bm{M}=\left(
\begin{array}{ccc}
 {\alpha_{1}} &  {\alpha_{2}}  &  {\alpha_{3}}\\
 {\alpha_{1}} &  {\alpha_{2}}  &  {\alpha_{3}}\\
 {\alpha_{1}} &  {\alpha_{2}}  &  {\alpha_{3}}
\end{array}
\right),
\eeqnn
with
\beqnn
\alpha_{k}=\Big(3V^*-2\lambda_{k}\Big)r_{k}.
\eeqnn
From \cite{HudSwe03}, the three eigenvalues of the matrix $\bm{A}$ are such that
for  a subcritical flow regime $(Fr<1)$,
\begin{align}
\lambda_{1} < 0  < \lambda_{2} \ll \lambda_{3}
\label{ineqlam}
\end{align}
and for a supercritical one $(Fr>1)$,
\begin{align}
 \lambda_{2} < 0  < \lambda_{1} < \lambda_{3}
\label{ineqlam2}
\end{align}
with $\lambda_{1}$ and $\lambda_{3}$ being the characteristic velocities of the water flow and $\lambda_{2}$ being
the characteristic velocity of the sediment motion. Obviously, the sediment motion is much slower than the
water flow.

 \subsection{Change of notations}
Hereafter, we consider the case where the flow regime is subcritical  and adopt the following notations:
$v(t,x)=\xi_1(t,x)$, $u_1(t,x)=\xi_2(t,x)$, $u_2(t,x)=\xi_3(t,x)$ and
coefficients  (characteristic velocities)
$\lambda_1=-\mu$, $\gamma_1=\lambda_2$ and $\gamma_2=\lambda_3$.
We introduce also the vector $\bm{u}=(u_{1},u_{2})^{tr}$,  the coefficients
 $\eta_{\rm{j}}= \alpha_{\rm{j+1}}$ for  $j=1,\,2$  and the matrix
\begin{align}
\bm{\sigma}=\left(
\begin{array}{cc}
 {\alpha_{2}}  &  {\alpha_{3}}\\
 {\alpha_{2}}  &  {\alpha_{3}}
\end{array}
\right).\label{sigma}
\end{align}

With the new variables, the set of equation (\ref{rimgen}) writes as:
\begin{subequations}
\begin{align}
& \partial_t u_{1}+\gamma_{1}\partial_x u_{1}=\sigma_{{11}}u_{1} +\sigma_{{12}}u_{\rm{2}} + \alpha_1 v \label{newvar1b} \\
& \partial_t u_{2}+\gamma_{2}\partial_x u_{2}=\sigma_{{21}}u_{1} +\sigma_{{22}}u_{\rm{2}} + \alpha_1 v \label{newvar1c}\\
& \partial_t v -\mu\partial_x v=\eta_{1}u_{1} + \eta_{2}u_{2} + \alpha_1 v. \label{newvar1a}
\end{align}
\label{newvar1}%
\end{subequations}
Introduce the variable $$w(t,x)=v(t,x)\exp\left(-\frac{\alpha_1}{\mu}x\right),$$
then the system (\ref{newvar1}) is transformed into
\begin{subequations}
\begin{align}
 \partial_t u_{1}+\gamma_{1}\partial_x u_{1}= & \sigma_{\rm{11}}u_{1} + \sigma_{\rm{12}}u_{2}\nonumber\\
&
+ \alpha_1 \exp\left(\frac{\alpha_1}{\mu}x\right)w\\
\partial_t u_{2}+\gamma_{2}\partial_x u_{2}=&  \sigma_{\rm{21}}u_{1} + \sigma_{\rm{22}}u_{2}\nonumber\\
&
+ \alpha_1 \exp\left(\frac{\alpha_1}{\mu}x\right)w\\
\partial_t w -\mu\partial_x w=&\eta_{1}\exp\left(\frac{\alpha_1}{\mu}x\right)u_{1}\nonumber\\
&+\eta_{2}\exp\left(\frac{\alpha_1}{\mu}x\right)u_{2}.
\end{align}
\label{newvar2}	
\end{subequations}
 We rewrite this system as:
\begin{subequations}
\begin{align}
& \partial_t u_{1}+\gamma_{1}\partial_x u_{1}=\sigma_{11}u_{1} +\sigma_{12}u_{2} + \alpha(x)w\\
& \partial_t u_{2}+\gamma_{2}\partial_x u_{2}=\sigma_{21}u_{1} +\sigma_{22}u_{2} + \alpha(x)w\\
& \partial_t w -\mu\partial_x w= \theta_{1}(x)u_{1}+\theta_{2}(x)u_{2} \label{newvar3_3}
\end{align}
\label{newvar3}%
\end{subequations}
with $\alpha(x)=\alpha_{1}\exp\left(\frac{\alpha_1}{\mu}x\right)$ and $\theta_{\rm{j}}(x) = \alpha_{\rm{j+1}}\exp\left(\frac{\alpha_1}{\mu}x\right)$ for $j=1,\,2$.

To close the writing of the system (\ref{newvar3}), we enclose to it the following boundary
and initial conditions
\begin{subequations}
\begin{align}
& u_{\rm{i}}(t,0)=q_{\rm{i}}w(t,0)\quad \text{for } {\rm{i}} =1,2, \label{bconda}\\
& w(t,1)=\rho_{1}u_{1}(t,1)+\rho_{2}u_{2}(t,1)+U(t), \label{bcondb}	\\
 & w(0,x)=w^{0}(x), \; u_{\rm{i}}(0,x)=u_{i}^{0}(x)\quad  \text{for } {\rm{i}} =1,2. \label{bcondc}	
\end{align}
\label{bcond}	
\end{subequations}
\begin{remm}
Let us mention that in the case where the flow regime is supercritical, the following changes
of variable will be considered (instead of the previous one)
$v(t,x)=\xi_2(t,x)$, $u_1(t,x)=\xi_1(t,x)$, $u_2(t,x)=\xi_3(t,x)$ and
coefficients
$\lambda_2=-\mu$, $\gamma_1=\lambda_1$ and $\gamma_2=\lambda_3$.
\end{remm}

\begin{figure}[H]
	\centering
	\includegraphics[width=3.7in]{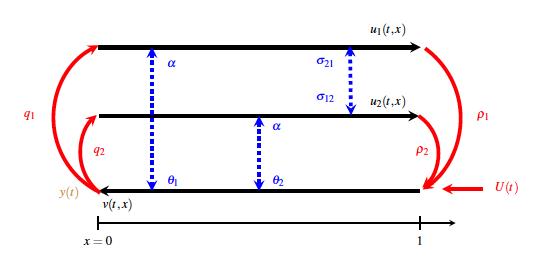}
	 \caption{ Schematic steep of the hyperbolic system. The internal coupling
	 between the states in the system and boundary conditions are depicted.}
	 \label{fig:originalsyst}
\end{figure}
 As in \cite{Dimeglio2013}, $u_1$, $u_2$ and $w$ are the distributed states and $U(t)$ is
 the control input as shown in Figure  \ref{fig:originalsyst}. The measured output is given by: $w(t,0)=y(t)$.
\section{Full State Controller Design}\label{backstepping}

\subsection{Backstepping transformation and target system}
Consider the following backstepping transformation
\begin{align}
&\psi_i (t,x)= u_{\rm{i}}(t,x)  \text{ for }  i=1,\,2\label{backtransfo2}\\
& \chi (t,x)= w(t,x)  - \int_{0}^{x} k_{1}(x,\xi)u_1(t,\xi)\,d\xi \nonumber  \\
&- \int_{0}^{x} k_{2}(x,\xi)u_2(t,\xi)\,d\xi
- \int_{0}^{x} k_{3}(x,\xi)w(t,\xi)\,d\xi.
\label{backtransfo1}
\end{align}
We now seek a sufficient condition on the functions $k_{\rm{i}}$ such that the transformation
(\ref{backtransfo2})-(\ref{backtransfo1}) maps the system (\ref{newvar3})-(\ref{bcond}) to the target
system
\begin{subequations}
\begin{align}
\partial_t \psi_1 + \gamma_{1}\partial_x \psi_1 =& \sigma_{11}\psi_1 + \sigma_{12}\psi_2 + \alpha(x) \chi \nonumber  \\
& + \int_{0}^{x} c_{11}(x,\xi)\psi_1(t,\xi)\,d\xi\nonumber\\
& + \int_{0}^{x}c_{12}(x,\xi)\psi_2(t,\xi)\,d\xi\nonumber  \\
&
+ \int_{0}^{x} \kappa_{1}(x,\xi)\chi(t,\xi)\,d\xi \label{targetsys_1}\\
 \partial_t \psi_2 + \gamma_{2}\partial_x \psi_2 =& \sigma_{21}\psi_1 + \sigma_{22}\psi_2 + \alpha(x) \chi\nonumber  \\
&+ \int_{0}^{x} c_{21}(x,\xi)\psi_1(t,\xi)\,d\xi \nonumber\\
&+ \int_{0}^{x}c_{22}(x,\xi)\psi_2(t,\xi)\,d\xi\nonumber  \\
&
+ \int_{0}^{x} \kappa_{2}(x,\xi)\chi(t,\xi)\,d\xi\label{targetsys_2}\\
 \partial_t  \chi- \mu \partial_x  \chi =&0
\end{align}
\label{targetsys}%
\end{subequations}
with the  following boundary conditions:
\begin{align}
\psi_{\rm{i}} (t,0) =q_{\rm{i}}\chi(t,0) \text{ for }  i=1,\,2   \text{ and  }  \chi(t,1)=0.
\label{bctarget}
\end{align}
This dynamic is schematically represented on Figure  \label{fig:targetsyst}.
In the system (\ref{targetsys}), $c_{\rm{ij}}(\cdot)$ and $\kappa_{\rm{i}}(\cdot)$ are functions
to be determined on the triangular domain
$$\mathbb{T}=\Big\{ (x,\xi)\in \mathbb{R}^2 | \; 0 \leq \xi  \leq x  \leq 1 \Big\}.$$
The system (\ref{targetsys})-(\ref{bctarget}) is designed as a copy of the original dynamics with the
coupling term in (\ref{newvar3_3}) removed. As will be shown later, the new terms in (\ref{targetsys_1}) and (\ref{targetsys_2}) are necessary for the design but they
will not affect the stability.
\begin{figure}[H]
	\centering
	\includegraphics[width=3.7in]{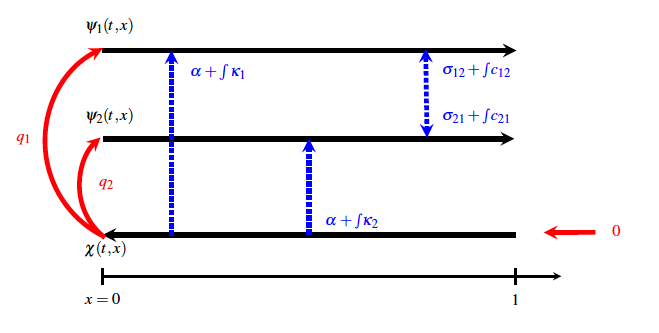}
	 \caption{Representation of the target system.}
	  \label{fig:targetsyst}
\end{figure}	
A sufficient condition for the transformation (\ref{backtransfo2})-(\ref{backtransfo1}) to map the original system 
(\ref{newvar3}) into the target system (\ref{targetsys}) is that the kernels
$k_{\rm{i}}$ satisfy the following system of first order hyperbolic PDEs:
\begin{subequations}
\begin{align}
&  \mu \partial_{x} k_1(x,\xi) - \gamma_1 \partial_{\xi} k_1(x,\xi) \nonumber  \\
&\quad=  \sigma_{11}k_{1}(x,\xi)
+ \sigma_{21}k_{2}(x,\xi) + \theta_{1}(\xi) k_3(x,\xi)\\
&  \mu \partial_{x} k_2(x,\xi) - \gamma_2 \partial_{\xi} k_2(x,\xi) \nonumber  \\
&\quad=  \sigma_{12}k_{1}(x,\xi)
+ \sigma_{22}k_{2}(x,\xi) + \theta_{2}(\xi) k_3(x,\xi)\\
& \mu  \partial_{x} k_3(x,\xi) + \mu \partial_{\xi} k_3(x,\xi) \nonumber  \\
&\quad=  \alpha(\xi)k_1(x,\xi)+ \alpha(\xi) k_2(x,\xi)
\end{align}
\label{kernelsyst}%
\end{subequations}
with the following boundary conditions:
\begin{subequations}
\begin{align}
&k_1(x,x)= -\frac{\theta_1(x)}{\gamma_1 + \mu},~
k_2(x,x)= - \frac{\theta_2(x)}{\gamma_2 + \mu},\\
&\mu k_3(x,0)= q_1\gamma_1k_1(x,0) + q_2\gamma_2k_2(x,0).
\end{align}
\label{kernelbc}	
\end{subequations}
  The existence, uniqueness and continuity of the solutions to the system (\ref{kernelsyst}) with boundary conditions (\ref{kernelbc})
  are assessed by Theorem 5.3 in \cite{Dimeglio2013}.

Besides, plugging (\ref{backtransfo2})-(\ref{backtransfo1}) into (\ref{targetsys}) and using  (\ref{newvar3})-(\ref{bcond}) yields for
$i=1,\,2$,
\begin{align}
 0 =& \int_{0}^{x}\Big [ \alpha(\xi) k_1(x,\xi) - c_{i1}(x,\xi)  \nonumber \\
&\qquad+ \int_{\xi}^{x} \kappa_1(x,s)k_{1}(s,\xi)\,ds \Big]u_1(\xi)\,d\xi\nonumber\\
& + \int_{0}^{x}\Big [ \alpha(\xi) k_2(x,\xi) - c_{i2}(x,\xi)  \nonumber \\
&\qquad+ \int_{\xi}^{x} \kappa_2(x,s)k_{2}(s,\xi)\,ds \Big]u_2(\xi)\,d\xi\\
& + \int_{0}^{x}\Big [ \alpha(\xi) k_3(x,\xi) - \kappa_{\rm i}(x,\xi)  \nonumber \\
&\qquad+ \int_{\xi}^{x} \kappa_{\rm i}(x,s)k_{3}(s,\xi)\,ds \Big]w(\xi)\,d\xi.
\label{coefkapc}	
\end{align}

The coefficients $\kappa_{\rm{i}}$ can be chosen to satisfy the following integral equation
for $i=1,\,2$
\begin{align}
 \kappa_{\rm{i}}(x,\xi)  = \alpha(x) k_3(x,\xi) + \int_{\xi}^{x} \kappa_{\rm i}(x,s)k_{3}(s,\xi)\,ds,
\label{kappa}	
\end{align}
and the coefficients $c_{\rm{ij}}$ can be chosen such that
\begin{align}
c_{\rm{ij}}(x,\xi) = \alpha(x) k_{\rm{j}}(x,\xi) &+ \int_{\xi}^{x}  \kappa_i(x,s)k_{\rm{j}}(s,\xi)\,ds \nonumber \\
&
 \text{for} \; i,\,j=1,\,2
\label{coefc}	
\end{align}
under the fact that the $k_{\rm i}$ exist and are  sufficiently smooth.

\subsection{Inverse transformation and control law}
To ensure that the target system and the  closed-loop system have equivalent stability properties,
the transformation (\ref{backtransfo2})-(\ref{backtransfo1}) has to be invertible. Since $\psi_{\rm i}=u_{\rm i}$, for $i=1,\,2$,
the transformation  (\ref{backtransfo1}) can be rewritten as
\begin{align}
 &\chi (t,x)+\int_{0}^{x} k_{1}(x,\xi)\psi_1(t,\xi)\,d\xi \nonumber\\
&+ \int_{0}^{x} k_{2}(x,\xi)\psi_2(t,\xi)\,d\xi\nonumber\\
&~~~~=
 w(t,x) - \int_{0}^{x} k_{3}(x,\xi)w(t,\xi)\,d\xi.
\label{invtransf}	
\end{align}
Let us define
\begin{align}
\Gamma (t,x) =& \chi (t,x)+\int_{0}^{x} k_{1}(x,\xi)\psi_1(t,\xi)\,d\xi \nonumber\\
&+ \int_{0}^{x} k_{2}(x,\xi)\psi_2(t,\xi)\,d\xi
\label{gammadef}	
\end{align}
Since $k_3$ is continuous  by Theorem $5.3$ in  \cite{Dimeglio2013}, there exists a unique continuous  inverse kernel
$l_3$ defined on $\mathbb{T}$, such that
\begin{align}
w(t,x) = \Gamma (t,x) + \int_{0}^{x} l_{3}(x,\xi)\Gamma(t,\xi)\,d\xi,
\label{inversegamma1}	
\end{align}
which yields the following  inverse transformation
Since $\psi_{\rm i}=u_{\rm i}$, for $i=1,\,2$, we could get the following relation from the first two equalities of $(\ref{newvar3})$ and $(\ref{targetsys})$:
\begin{align}
&\alpha(x)w= \alpha(x) \chi  + \int_{0}^{x} c_{11}(x,\xi)\psi_1(t,\xi)\,d\xi\nonumber  \\
& + \int_{0}^{x}c_{12}(x,\xi)\psi_2(t,\xi)\,d\xi
+ \int_{0}^{x} \kappa_{1}(x,\xi)\chi(t,\xi)\,d\xi.
\end{align}
Thus, we could write the following
\begin{align}
&w(t,x) =  \chi (t,x)+\int_{0}^{x} l_{1}(x,\xi)\psi_1(t,\xi)\,d\xi \nonumber  \\
&+\int_{0}^{x} l_{2}(x,\xi)\psi_2(t,\xi)\,d\xi
+ \int_{0}^{x} l_{3}(x,\xi)\chi(t,\xi)\,d\xi,
\label{inversegamma2}	
\end{align}
where for $i=1,\,2$,
\begin{align}
l_{\rm i}(x,\xi)=  k_{\rm i} (t,\xi)+\int_{\xi}^{x} k_{\rm i} (x,\xi)l_3(\xi,s)\,ds.
\label{ldef}	
\end{align}
Thus, the control  law $U(t)$ can be obtained by plugging the transformation  (\ref{backtransfo1}) into
(\ref{newvar3}). Readily, $\chi(t,1)=0$ implies that
\begin{align}
U(t)=& -\rho_{1}u_{1}(t,1) - \rho_{2}u_{2}(t,1) \nonumber\\
&+ \int_{0}^{1}\Big [ k_1(1,\xi)u_{1}(x,\xi) \nonumber \\
& +  k_2(1,\xi)u_{2}(x,\xi) +k_{3}(1,\xi)w(1,\xi) \Big]\,d\xi.
\label{controldef}
\end{align}
The $k_i$ in the integral term designate the kernel functions and satisfy the system  (\ref{kernelsyst})-(\ref{kernelbc}).
\subsection{ Stability of the target system and the closed-loop control system}

We first prove exponential stability of the target system
(\ref{targetsys})-(\ref{bctarget}).
\begin{lemm}
For any given initial condition
$(\psi^0_1, \ \psi^0_2, \ \chi^0)^T \, \in \, \left(\mathcal{L}^{2}([0,1])\right)^3$ and under the assumption that
$c_{\rm{ij}},~\kappa_{\rm{i}} \, \in \mathcal{C}(\mathbb{T})$, the equilibrium $(\psi_1, \ \psi_2, \ \chi)^T=(0, \ 0, \ 0)^T$ of the target system
(\ref{targetsys})-(\ref{bctarget}) is  $\mathcal{L}^2$-exponentially stable.
\end{lemm}
\begin{pf}
 The stability proof  is  based on  the time  differentiation of  the following Lyapunov function:
  \begin{align}
V_1(t)=& \int_0^1 a_1 e^{-\delta_1 x}\left(\frac{\psi^2_{1}(t,x)}{\gamma_1} + \frac {\psi^2_{2}(t,x)}{\gamma_2}\right)dx\nonumber  \\
 &+ \int_0^1  \frac{1+x  }{\mu}\chi^2(t,x) dx, \label{V1}
\end{align}
where $a_1$ and $\delta_1$ are strictly positive parameters
to be determined.

Differentiating this function with respect to time, we get:
\begin{align}
\dot V_1(t)=&2 \int_0^1 a_1 e^{-\delta_1 x}\left(\frac{  \psi_{1}  \partial_t  \psi_{1}}{\gamma_1}+ \frac { \psi_{2}  \partial_t  \psi_{2}}{\gamma_2}\right)dx \nonumber  \\
&+ 2\int_0^1  \frac{1+x}{\mu}  \chi \partial_t  \chi dx.\label{dotV1}
\end{align}
By taking into account the target system (\ref{targetsys})-(\ref{bctarget}) and integrating by parts, we have
 \begin{align}
\dot V_1(t)=& \Big[ -a_1 e^{-\delta_1 x} (\psi^2_{1}(t,x)  + \psi^2_{2}(t,x))\nonumber\\
&~~~~~~+ (1+x) \chi^2(t,x)\Big]_0^1 \nonumber  \\
&-\int_0^1 \chi^2(t,x)\ dx+ \int_0^1  a_1 e^{-\delta_1 x}  \nonumber\\
& \times \Psi^{T}(t,x) \left( -\delta_1 I_2 + 2\bm{\Gamma}_{inv}\bm{\sigma} \right) \Psi (t,x) \  dx\nonumber\\
 &+ 2 \int_0^1 a_1 e^{-\delta_1 x} \Psi^{T}(t,x) \bm{\Gamma}_{inv} \bm{\alpha}(x) \chi (t,x)  \ dx \nonumber \\
&+  2\int_0^1   a_1 e^{-\delta_1 x}  \int_{0}^{x} \Psi^{T}(t,x) \bm{\Gamma}_{inv} \nonumber \\
&\times \left(\bm{C}(x, \xi) \Psi (t,\xi)\right. \left. +K(x, \xi) \chi(t,\xi) \right) \ d\xi \ dx\nonumber,
\end{align}
where the matrix $\bm{\sigma}$ is defined in $(\ref{sigma})$, the vectors $\Psi(t,x)$, $\bm{\alpha}(x)$, $K(x,\xi)$ and the matrices $\bm{\Gamma}_{inv}$, $\bm{C}(x, \xi)$ are given by
\begin{align}
&{\Psi}(t,x)=\left(
\begin{matrix}
 { \psi_1(t,x)}  \\
 { \psi_2(t,x)}
\end{matrix}
\right),~
\bm{\alpha}(x)=\left(
\begin{matrix}
 \alpha(x)  \\
 \alpha(x)
\end{matrix}
\right)
\label{psialphaK}\\
&K(x,\xi)=\left(
\begin{matrix}
 { \kappa_1(x,\xi)}  \\
 { \kappa_2(x,\xi)}
\end{matrix}
\right),
\bm{\Gamma}_{inv}=\left(\begin{matrix}
                                       \frac{1}{\gamma_1} & 0 \\
                                       0 & \frac{1}{\gamma_2}
                                     \end{matrix}
                                   \right)\\
&\bm{C}(x, \xi)=\left(
\begin{matrix}
 { c_{11}}(x, \xi) &  { c_{12}}(x, \xi) \\
 { c_{21}}(x, \xi) &  { c_{22}}(x, \xi)
\end{matrix}
\right).\label{gammaC}
\end{align}
Assume that for $M>0$ and $\epsilon>0$, we have
\begin{align}
& \|\bm{\sigma}  \|, \ \|\bm{\alpha}(x)\|,  \|\bm{C}(x,\xi) \|,  \  \|K(x,\xi) \|\leq M,\nonumber\\
 &~~~~~~~~~~~~~~~~~~~~\forall x\in [0,1], \xi \in [0,x],\\
& \gamma_i(x)>\epsilon, \forall  i=1,2, \forall x\in [0,1],
\end{align}
where the matrix/vector norms $\|\cdot\|$ are compatible with the other corresponding matrix/vector norms. 
Hence, using Young's inequalities the following relations are derived
 \begin{align}
\hspace{-2cm}
& 2\int_0^1  a_1 e^{-\delta_1 x}  \Psi^{T}(t,x)  \bm{\Gamma}_{inv}\bm{\sigma}\Psi (t,x) \  dx\nonumber\\
 &~~~~\leq 2\frac{M}{\epsilon} \int_0^1 a_1 e^{-\delta_1 x} \Psi^{T}(t,x)  \Psi (t,x) \  dx
\end{align}
 \begin{align}
& 2 \int_0^1 a_1 e^{-\delta_1 x} \Psi^{T}(t,x) \bm{\Gamma}_{inv} \bm{\alpha}(x) \chi (t,x)  \ dx \nonumber\\
&~~~~\leq \int_0^1 a_1 e^{-\delta_1 x} \left(\Psi^{T}(t,x) \bm{\Gamma}_{inv} \bm{\alpha}(x)\right.\nonumber\\
&\left.~~~~~~~~~~~~~~\times \bm{\alpha}^T(x) \bm{\Gamma}_{inv}\Psi(t,x) +\chi^2 (t,x)\right) dx \nonumber \\
&~~~~\leq a_1 \left(\frac{M}{\epsilon}\right)^2 \int_0^1 e^{-\delta_1 x} \Psi^{T}(t,x) \Psi(t,x) dx\nonumber\\
&~~~~~~~~+a_1 \int_0^1 e^{-\delta_1 x} \chi^2 (t,x) dx
\end{align}
and
\begin{align}
&2\int_0^1   a_1 e^{-\delta_1 x}  \int_{0}^{x} \Psi^{T}(t,x) \bm{\Gamma}_{inv}\bm{C}(x, \xi) \Psi (t,\xi)  \ d\xi \ dx \nonumber\\
&\leq \frac{M}{\epsilon} \int_0^1 a_1 e^{-\delta_1 x} \int_0^x \left(\Psi^{T}(t,x)\Psi(t,x)\right.\nonumber\\
&~~~~\left.+\Psi^{T}(t,\xi)\Psi(t,\xi)\right) d\xi dx\nonumber\\
&= \frac{M}{\epsilon} \int_0^1 a_1 e^{-\delta_1 x} x \Psi^{T}(t,x)\Psi(t,x) dx\nonumber\\
&~~~~+ \frac{M}{\delta_1 \epsilon} \int_0^1 a_1 \left(e^{-\delta_1 x}-e^{-\delta_1 }\right) \Psi^{T}(t,x)\Psi(t,x) dx\nonumber\\
&\leq a_1\int_0^1 e^{-\delta_1 x}\left(\frac{M}{\epsilon}x+\frac{M}{\delta_1 \epsilon} \right) \Psi^{T}(t,x)\Psi(t,x) dx\end{align}
\begin{align}
&2\int_0^1   a_1  e^{-\delta_1 x}  \int_{0}^{x} \Psi^{T}(t,x) \bm{\Gamma}_{inv} K(x, \xi) \chi(t,\xi)  \ d\xi \ dx\nonumber\\
&\leq a_1\int_0^1 e^{-\delta_1 x}\int_0^x \left(\Psi^{T}(t,x) \bm{\Gamma}_{inv} K(x, \xi)\right.\nonumber\\
&~~~~\left.\times K^T(x, \xi)\bm{\Gamma}_{inv}\Psi(t,x) + \chi^2(t,\xi)\right)d\xi dx\nonumber\\
&\leq a_1 \left(\frac{M}{\epsilon}\right)^2 \int_0^1 e^{-\delta_1 x}x\Psi^{T}(t,x) \Psi(t,x) dx\nonumber\\
&~~~~ + a_1 \frac{1}{\delta_1} \int_0^1 e^{-\delta_1 x} \chi^2(t,x) dx.
\end{align}
Thus, using the boundary conditions (\ref{bctarget}), we obtain the following inequality
 \begin{align}
\dot V_1(t)\leq& \left( a_1 \sum_{i=1}^2 q_i^2 -1\right) \chi^2(t,0) \nonumber  \\
&-\int_0^1 \left(1-a_1(1+\frac{1}{\delta_1})e^{-\delta_1 x}\right)\chi^2(t,x)\ dx  \nonumber\\
&- a_1 \int_0^1 e^{-\delta_1 x} \Psi^{T}(t,x) P_1(x)\Psi (t,x) \  dx,
\end{align}
where
\begin{align}
P(x)=&\left(\delta_1-2\frac{M}{\epsilon}-\frac{M}{\epsilon}x-2\left(\frac{M}{\epsilon}\right)^2-\frac{M}{\delta_1 \epsilon}\right)I_2\nonumber  \\
&-2\bm{\Gamma}_{inv}\bm{\sigma}.
\end{align}

First, we choose the tuning parameter $\delta_1>0$  sufficiently large so that the matrix $P(x), x\in [0,1]$ is positive definite. Then, by choosing
\begin{align}
0< a_1 < \min\left\{   \frac{1}{\sum\limits_{i=1}^2 q_i^2}, \frac{\delta_1}{\delta_1+1}   \right\},
\end{align}
we could derive exponential stability of the target system.
\end{pf}

Then, from the continuity and invertibility of the backstepping transformation (\ref{backtransfo2})-(\ref{backtransfo1}), 
we could derive equivalence between the original system (\ref{newvar3}) (with the boundary and initial conditions (\ref{bcond}) and the control law
(\ref{controldef})) and the target system (\ref{targetsys})-(\ref{bctarget}). Thus, the following theorem is proved.
\begin{thm}
 Consider the system (\ref{newvar3}) with the boundary and initial conditions (\ref{bcond}) and the control law
(\ref{controldef}). Then under the assumptions that  the initial data are in $\left(\mathcal{L}^2([0,1])\right)^3$, the origin is exponentially
stable in the $\mathcal{L}^2$ sense.
\end{thm}

\section{Backstepping Observer Design}\label{observer}
The feedback  controller  (\ref{controldef}) requires a full state measurement  across the spatial domain.
In this section we are interested in the design of  a boundary state observer for estimation of  the distributed  states of the system
(\ref{newvar3})-(\ref{bcond}) over the whole spatial domain using the measured output $w(t,0)=y(t)$.  The observer

\begin{subequations}
\begin{align}
& \partial_t \hat u_{1}+\gamma_{1}\partial_x \hat u_{1}=\sigma_{11}\hat u_{1} +\sigma_{12}\hat u_{2} + \alpha(x)\hat w \nonumber\\
&~~~~~~~~~~~~~~~~~~~ -p_1(x)[y(t)-\hat w(t,0)]\\
& \partial_t \hat u_{2}+\gamma_{2}\partial_x \hat  u_{2}=\sigma_{21} \hat u_{1} +\sigma_{22} \hat u_{2} + \alpha(x)\hat w\nonumber\\
&~~~~~~~~~~~~~~~~~~~-p_2(x)[y(t)-\hat w(t,0)]\\
& \partial_t \hat w -\mu\partial_x\hat  w= \theta_{1}(x)\hat u_{1}+\theta_{2}(x)\hat u_{2} \nonumber\\
&~~~~~~~~~~~~~~~~~~~-p_3(x)[y(t)-\hat w(t,0)],
\end{align}
\label{newvar3-obs}%
\end{subequations}
\noindent where $(  \hat u_1, \ \hat u_2, \  \hat w)^T$  is  the estimated state vector, consists of a copy of the plant plus an output injection and mimics the well-known finite dimensional observer format. The functions $\theta_{\rm{j}}(x) = \alpha_{\rm{j+1}}$ for $j=1,\,2$ and  $\alpha(x)$ are the ones defined for the transformed  system (\ref{newvar3}).
The following boundary  conditions  have to be  considered:
\begin{subequations}
\begin{align}
& \hat u_{\rm{i}}(t,0)=q_{\rm{i}}y(t)\quad \text{for } {\rm{i}} =1,2 \label{bconda-obs}\\
& \hat w(t,1)=\rho_{1}\hat u_{1}(t,1)+\rho_{2}\hat u_{2}(t,1)+U(t). \label{bcondb-obs}	
\end{align}
\label{bcond-obs}	
\end{subequations}
Our objective is to find  $p_1(x)$, $p_2(x)$ and $p_3(x)$ such that the estimated state vector $(\hat w, \  \hat u_1, \ \hat u_2 )$
converges to the real state vector $( w, \   u_1, \ u_2 )$ in finite time. Defining
\begin{align}
&\left(
  \begin{array}{ccc}
    \tilde w & \tilde u_1 & \tilde u_2 \\
  \end{array}
\right)^T
= \left(
  \begin{array}{ccc}
    w-\hat w & u_1-\hat u_1 &  u_2-\hat u_2 \\
  \end{array}
\right)^T\label{error}
 \end{align}
as the error variable vector, we obtain the following error system
\begin{subequations}
\begin{align}
 \partial_t \tilde w -\mu\partial_x\tilde  w=& \theta_{1}(x)\tilde u_{1}+\theta_{2}(x)\tilde u_{2} +p_3(x)\tilde w(t,0)\\
 \partial_t \tilde u_{1}+\gamma_{1}\partial_x\tilde u_{1}=&\sigma_{11}\tilde u_{1} +\sigma_{12}\tilde u_{2} + \alpha(x)\tilde w\nonumber  \\
& +p_1(x)\tilde w(t,0)\\
 \partial_t \tilde u_{2}+\gamma_{2}\partial_x\tilde  u_{2}=&\sigma_{21} \tilde u_{1} +\sigma_{22} \tilde u_{2} + \alpha(x)\tilde w\nonumber  \\
&+p_2(x)\tilde w(t,0)
\end{align}
\label{newvar3-obs-error}%
\end{subequations}
with the boundary conditions
\begin{subequations}
\begin{align}
& \tilde w(t,1)=\rho_{1}\tilde u_{1}(t,1)+\rho_{2}\tilde u_{2}(t,1), \label{bconda-obs-error}\\
& \tilde u_{\rm{i}}(t,0)=0\quad \text{for } {\rm{i}} =1,2. \label{bcondb-obs-error}
\end{align}	
\label{bcond-obs-error}
\end{subequations}
\subsection{Backstepping transformation and the target error system}
Similarly to the controller design, we use the following invertible backstepping transformation	
\begin{subequations}
\begin{align}
& \tilde u_{\rm i} (t,x)= \tilde \pi_{\rm i}(t,x)  + \int_{0}^{x} m_{\rm i}(x,\xi)\tilde \phi(t,\xi)\,d\xi\label{obs-1-trans-bis}\\
&~~~~~~~~~~~~~~~~~~~~~~~~~~~~~ \quad \rm{for } ~{\rm{i}} =1,2\nonumber\\
&  \tilde w (t,x)= \tilde \phi(t,x)  +\int_{0}^{x} m_{\rm 3}(x,\xi)\tilde \phi(t,\xi)\,d\xi, \label{obs-1-trans}
 \end{align}
\label{backtransfo-obs}%
\end{subequations}
 where the kernels $m_{\rm{i}} (\cdot) \  \text{ for }  i=1,\,2,\,3 $ are defined in the triangular domain $\mathbb{T}$  to map the
error system  (\ref{newvar3-obs-error})-(\ref{bcond-obs-error}) into the following exponentially stable target system
 \begin{subequations}
\begin{align}
&\partial_t \tilde \pi_1 + \gamma_{1}\partial_x \tilde \pi_1 = \sigma_{11}\tilde \pi_1 + \sigma_{12}\tilde \pi_2\nonumber\\
&+ \int_{0}^{x} g_{11}(x,\xi)\tilde \pi_1(t,\xi)\,d\xi + \int_{0}^{x}g_{12}(x,\xi)\tilde \pi_2(t,\xi)\,d\xi\\
&
\partial_t \tilde \pi_2 + \gamma_{2}\partial_x \tilde \pi_2 = \sigma_{21} \tilde \pi_1 + \sigma_{22}\tilde \pi_2\nonumber\\
&+ \int_{0}^{x} g_{21}(x,\xi)\tilde \pi_1(t,\xi)\,d\xi + \int_{0}^{x}g_{22}(x,\xi)\tilde \pi_2(t,\xi)\,d\xi\\
& \partial_t \tilde \phi - \mu \partial_x \tilde \phi = \theta_{1}(x)\tilde \pi_{1}+\theta_{2}(x)\tilde \pi_{2}\nonumber\\
&+ \int_{0}^{x}h_{1}(x,\xi)\tilde \pi_1(t,\xi)\,d\xi  +\int_{0}^{x}h_{2}(x,\xi)\tilde \pi_2(t,\xi)\,d\xi
\end{align}
\label{targetsys-obs}%
\end{subequations}
with the  boundary conditions
 \begin{subequations}
\begin{align}
 &\tilde \pi_{\rm{i}} (t,0) =0 \text{ for }  i=1,\,2\\
  &\tilde \phi (t,1)=\rho_{1} \tilde \pi_{1}(t,1)+\rho_{2}\tilde \pi_{2}(t,1).
\end{align}\label{bctarget-obs}
\end{subequations}
Here the functions $g_{ij}$  and  $h_i$ have to be determined on the triangular domain $\mathcal{T}$.
As previously, we are attempting to find some sufficient condition for the kernels to match the target system.
Differentiating the transformations  (\ref{backtransfo-obs})  in time and space and  substituting  the results into
(\ref{newvar3-obs-error}) with the help of (\ref{targetsys-obs}), the following PDEs are derived for the kernels
\begin{subequations}
\begin{align}
\gamma_1 \partial_x m_{1}-\mu \partial_\xi m_{1} &= \sigma_{11} m_{1} +\sigma_{12} m_{2} + \alpha(x)m_3 ,\\
\gamma_2  \partial_x m_{2}-\mu \partial_\xi m_{2} &= \sigma_{21} m_{1} +\sigma_{22} m_{2} + \alpha(x)m_3,\\
\mu  \partial_x m_{3} + \mu \partial_\xi m_3 &= -\theta_{1}(x) m_{1}-\theta_{2}(x) m_{2}.
\end{align}
\label{obs-error-kernel}	
\end{subequations}
To close the writing of the  above system, the following boundary conditions are imposed:
\begin{subequations}
\begin{align}
& m_1(x,x)=\frac{1}{\gamma_1+ \mu} \alpha(x)\\
&m_2(x,x)=\frac{1}{\gamma_2+ \mu} \alpha(x)\\
& m_3(1,\xi)=\rho_1 m_1(1,\xi) +\rho_2 m_2(1,\xi).
\end{align}
\label{bc-obs-error-kernel}
\end{subequations}
The observer gains  are  defined by
\begin{align}
 p_i(x)=\mu m_i(x,0)  \text{ for }  i=1,\,2, \ 3,
\label{gain-obs}
\end{align}
and the integral coupling coefficients are defined by the following equations:\\
 \begin{subequations}
\begin{align}
 &h_i(x,\xi)=-\theta(\xi) m_3(x,\xi)\nonumber\\
&- \int_{\xi}^{x}m_{3}(x,s)h_i(s,\xi )\,ds, ~ \rm{ for }~ i =1,\,2, \\
&g_{i,j}(x,\xi)=-\theta_j (\xi) m_i(x,\xi)\nonumber\\
&-\int_{\xi}^{x}m_{i}(x,s)h_j(s,\xi )\,ds,~ \rm{ for } ~\{ i,\ j\} =1,\,2.
\end{align}
\label{obs-int-coupling}	
\end{subequations}

\subsection{Inverse Transformation}
The  continuity of the kernel $m_3$ in the transformation (\ref{obs-1-trans}) guarantees the existence
of a unique continuous  inverse kernel $r_3$ in the transformation
\begin{align}
 \tilde \phi (t,x)= \tilde w(t,x)  +\int_{0}^{x} r_{\rm 3}(x,\xi)\tilde w(t,\xi)\,d\xi \label{obs-1-trans-inv}
 \end{align}
define on $\mathbb{T}$, and
\begin{align}
 r_3(x,\xi)= -m_3 (x,\xi)  -\int_{\xi}^{x} m_{3}(x,s)r_3(s, \xi)\,ds. \label{obs-1-trans-inv-kernel}
 \end{align}
Substituting  (\ref{obs-1-trans-inv-kernel})  into (\ref{obs-1-trans-bis}),  we obtain
\begin{align}
 \tilde\pi_{\rm i}(t,x)=&  \tilde u_{\rm i}(t,x) -\int_{0}^{x} m_{\rm i}(x,\xi)\tilde w(t,\xi)\,d\xi\nonumber\\
&- \int_{0}^{x}\int_{0}^{\xi}m_{\rm i}(x,\xi) r_3( \xi,s)\tilde w(t,s) ds d\xi \nonumber\\
=&  \tilde u_{\rm i}(t,x) -\int_{0}^{x} m_{\rm i}(x,\xi)\tilde w(t,\xi)\,d\xi\nonumber\\
&- \int_{0}^{x} \tilde w(t,\xi)  \int_{\xi}^{x}m_{\rm i}(x,s) r_3(s, \xi) \ ds d\xi,\nonumber\\
&~~~~~~~~~~~~~~~~~~~~~~~~~~~~~ \text{ for }  i =1,\,2,\nonumber
 \end{align}
and hence for $i =1,\,2, $
\begin{align}
 \tilde\pi_{\rm i}(t,x)&=  \tilde u_{\rm i}(t,x) +\int_{0}^{x} r_{\rm i}(x,\xi)\tilde w(t,\xi)\,d\xi
 \end{align}
where
\begin{align}
 r_{\rm i}(x,\xi)&=  -\tilde m_{\rm i}(x,\xi) -\int_{\xi}^{x} m_{\rm i}(x,s)r_3(s,\xi)\,ds. \nonumber
 \end{align}

\subsection{Stability of the target error system and convergence of the designed observer}
The observer is exponentially convergent to the original system. We first prove exponential stability of the target error system (\ref{targetsys-obs}).
\begin{lemm}
Under the assumptions $\psi^0_1$,  $\psi^0_2$,   $\chi^0 \in \mathcal{L}^2([0,1])$  and  $g_{\rm{ij}}$,
$h_{\rm{i}} \in \ \mathcal{C}(\mathbb{T})$,
the  system (\ref{targetsys-obs}) with
boundary conditions  (\ref{bctarget-obs}) and given initial condition $(\psi^0_1,  \psi^0_2,  \chi^0)$ is exponentially
stable in the $\mathcal{L}^2$ sense.
\end{lemm}
\begin{pf}
 The stability proof  is  based on  the time  differentiation of  the following Lyapunov function
\begin{align}
V_2(t)&= \int_0^1 a_2 e^{-\delta_2 x}\left(\frac{\tilde \pi^2_{1}(t,x)}{\gamma_1}
+ \frac {\tilde \pi^2_{2}(t,x)}{\gamma_2}\right)dx \nonumber  \\
&+ \int_0^1  \frac{e^{\delta_2 x}  }{\mu} \tilde\phi^2(t,x) dx, \label{V2}
\end{align}
where $a_2$ and $\delta_2$ are strictly positive parameters
to be determined.
Differentiating this function with respect to time, we get:
\begin{align}
\dot V_2(t)&=2 \int_0^1 a_2 e^{-\delta_2 x}\left(\frac{ \tilde \pi_{1}  \partial_t \tilde \pi_{1}}{\gamma_1}
+ \frac {\tilde \pi_{2}  \partial_t \tilde \pi_{2}}{\gamma_2}\right)dx\nonumber  \\
& + 2\int_0^1  \frac{e^{\delta_2 x}}{\mu} \tilde \phi \partial_t \tilde \phi dx.\label{dotV2}
\end{align}
Taking into account of the target system (\ref{targetsys-obs}) and integrating by parts, we rewrite (\ref{dotV2}) as
\begin{align}
&\dot V_2(t)= \Big[ -a_2 e^{-\delta_2 x} (\tilde\pi^2_{1}(t,x)  + \tilde\pi^2_{2}(t,x))\nonumber  \\
&~~~~~~~~+ e^{\delta_2 x} \tilde\phi^2(t,x)\Big]_0^1- \delta_2  \int_0^1  e^{\delta_2 x} \tilde \phi^2(t,x)\ dx\nonumber\\
&+ 2\int_0^1  a_2 e^{-\delta_2 x}  \Pi^{T}(t,x){\bm \Gamma}_{inv} {\bm \sigma} \   \Pi (t,x) dx\nonumber  \\
&-\delta_2  \int_0^1  a_2 e^{-\delta_2 x}\Pi^{T}(t,x)  \   \Pi (t,x) \ dx\nonumber \\
 &+ 2 \int_0^1\int_{0}^{x}  a_2 e^{-\delta_2 x} \Pi^{T}(t,x) {\bm \Gamma}_{inv}\nonumber\\
&~~~~~~~~ \times{\bf G}(x, \xi) \   \Pi (t,\xi) d\xi \ dx\nonumber  \\
&+ 2\int_0^1  \frac{e^{\delta_2 x}}{\mu}  \tilde\phi(t,x) {\bm \theta}(x) \Pi(t,x) dx \nonumber \\
&+ 2\int_0^1  \frac{e^{\delta_2 x}}{\mu} \tilde \phi(t,x) \int_{0}^{x}  \bm{h}(x, \xi) \Pi(t,\xi) \,d\xi  \ dx \nonumber,
\end{align}
where the matrices ${\bm \Gamma}_{inv}, {\bm \sigma}$ are defined by $(\ref{gammaC})$ and $(\ref{sigma})$,
the vector $\Pi, {\bm \theta, \bm h}$ and the matrix $\bm{G}$ are given by
\begin{align}
&{\Pi}(t,x)=\left(
\begin{matrix}
 { \tilde\pi_1(t,x)}  \\
 { \tilde\pi_2(t,x)}
\end{matrix}
\right),\\
&\bm{G}(x, \xi)=\left(
\begin{matrix}{cc}
 g_{11}(x, \xi) &  g_{12}(x, \xi) \\
 g_{21}(x, \xi) &  g_{22}(x, \xi)
\end{matrix}
\right),\\
&{\bm \theta}(x)=\left(
\begin{matrix}{cc}
 { \theta_1(x))}  &  { \theta_2(x)}
\end{matrix}
\right),\\
&\bm{h}(x, \xi)=\left(
\begin{array}{cc}
 \bm{h}_{1}(x, \xi) &  \bm{h}_{2}(x, \xi)
\end{array}
\right).
\end{align}
Assume that for $\tilde M>0$, we have
\begin{align}
& \|\bm{G}(x, \xi) \|, \ \|\bm{\theta}(x)\|,  \|\bm{h}(x,\xi) \| \leq \tilde M,\nonumber\\
&~~~~~~~~~~~~~~~~~~~~ \forall x\in [0,1], \xi \in [0,x],
\end{align}
where the matrix/vector norms $\|\cdot\|$ are compatible with the other corresponding matrix/vector norms. Hence, using Young's inequality, the following are derived
 \begin{align}
\hspace{-2cm}
& 2\int_0^1  a_2 e^{-\delta_2 x}  \Pi^{T}(t,x)  \bm{\Gamma}_{inv}\bm{\sigma}\Pi (t,x) \  dx\nonumber\\
 &~~~~\leq 2\frac{\tilde M}{\epsilon} \int_0^1 a_2 e^{-\delta_2 x} \Pi^{T}(t,x)  \Pi (t,x) \  dx
\end{align}
and,
\begin{align}
& 2\int_0^1  \frac{e^{\delta_2 x}}{\mu}  \tilde\phi(t,x) {\bm \theta}(x) \Pi(t,x) dx\nonumber\\
&~~~~\leq \int_0^1 \frac{e^{\delta_2 x}}{\mu} \left(\tilde\phi^2 (t,x)\right.\nonumber\\
&\left.~~~~~~~~~~+\Pi^{T}(t,x) {\bm \theta}^T(x) {\bm \theta}(x)\Pi(t,x) \right) dx \nonumber\\
&~~~~\leq \int_0^1 \frac{e^{\delta_2 x}}{\mu} \tilde\phi^2 (t,x) dx\nonumber\\
&~~~~~~~~~~+\frac{\tilde M^2}{\mu}\int_0^1 e^{\delta_2 x} \Pi^{T}(t,x) \Pi(t,x) dx
\end{align}
and,
\begin{align}
&2 \int_0^1\int_{0}^{x}  a_2 e^{-\delta_2 x} \Pi^{T}(t,x) {\bm \Gamma}_{inv}{\bf G}(x, \xi) \   \Pi (t,\xi) d\xi \ dx \nonumber\\
&~~~~\leq \frac{\tilde M}{\epsilon} \int_0^1 a_2 e^{-\delta_2 x} \int_0^x \left(\Pi^{T}(t,x)\Pi(t,x)\right.\nonumber\\
&\left.~~~~~~~~~+\Pi^{T}(t,\xi)\Pi(t,\xi)\right) d\xi dx\nonumber\\
&~~~~= \frac{\tilde M}{\epsilon} \int_0^1 a_2 e^{-\delta_2 x} x \Pi^{T}(t,x)\Pi(t,x) dx\nonumber\\
&~~~~~~+ \frac{\tilde M}{\delta_2 \epsilon} \int_0^1 a_2 \left(e^{-\delta_2 x}-e^{-\delta_2 }\right) \Pi^{T}(t,x)\Pi(t,x) dx\nonumber\\
&~~~~\leq a_2\int_0^1 e^{-\delta_2 x}\left(\frac{\tilde M}{\epsilon}x+\frac{\tilde M}{\delta_2 \epsilon} \right) \Pi^{T}(t,x)\Pi(t,x) dx
\end{align}
and finally
\begin{align}
&2\int_0^1  \frac{e^{\delta_2 x}}{\mu} \tilde \phi(t,x) \int_{0}^{x}  \bm{h}(x, \xi) \Pi(t,\xi) \,d\xi  \ dx\nonumber\\
&~~~~\leq \int_0^1 \frac{e^{\delta_2 x}}{\mu} \int_0^x \left(\Pi^{T}(t,\xi) \bm{h}^T(x, \xi)\bm{h}(x, \xi)\Pi(t,\xi)\right.\nonumber\\
&\left.~~~~~~~~~~~~~~~~ + \tilde\phi^2(t,x)\right)d\xi dx\nonumber\\
&~~~~\leq  \frac{\tilde M^2}{\delta_2 \mu} \int_0^1 \left( e^{\delta_2}-e^{\delta_2 x}\right) \Pi^{T}(t,x) \Pi(t,x) dx \nonumber\\
&~~~~~~~~ +\frac{1}{\mu} \int_0^1 e^{\delta_2 x} x \tilde\phi^2(t,x) dx.
\end{align}
Thus,
 \begin{align}
\hspace{-1cm}
\dot V_2(t)\leq& - e^{-\delta_2 }a_2\left(\tilde\pi^2_{1}(t,1)  + \tilde\pi^2_{2}(t,1) \right)\nonumber\\
&+e^{\delta_2 }\left( \rho_1\tilde\pi_{1}(t,1)  +\rho_2 \tilde\pi_{2}(t,1) \right)^2\nonumber\\
&- \int_0^1   e^{\delta_2 x}\left(\delta_2-\frac{1+x}{\mu}\right)\tilde \phi^2(t,x) \ dx\nonumber\\
&- \int_0^1\Pi^{T}(t,x)\left[ e^{-\delta_2 x}  \left(\delta_2 a_2\right.\right.\nonumber\\
&\left.\left.-{\frac{a_2(2+x+1/\delta_2)\tilde M}{\epsilon}} \right)- \frac{\tilde M^2}{\delta_2 \mu}   e^{\delta_2}\right.\nonumber\\
&\left.+e^{\delta_2 x}\left(\frac{1}{\delta_2}-1\right)\frac{\tilde M^2}{\mu}\right]\   \Pi (t,x) dx.\nonumber
\end{align}
With the help of the boundary conditions (\ref{bctarget-obs}), we obtain
\begin{align}
\dot V_2(t)
\leq& - e^{-\delta_2 }\left[\left(a_2-2\rho_1^2e^{2\delta_2}\right)\tilde\pi^2_{1}(t,1) \right.\nonumber\\
&\left. + \left(a_2-2\rho_2^2e^{2\delta_2}\right)\tilde\pi^2_{2}(t,1)\right]\nonumber\\
&- \int_0^1   e^{\delta_2 x}\left(\delta_2-\frac{1+x}{\mu}\right)\tilde \phi^2(t,x) \ dx\nonumber\\
&-\int_0^1\Pi^{T}(t,x)e^{-\delta_2 x}\tilde P(x)\   \Pi (t,x) dx,
\label{V2dot}
\end{align}
where
\begin{align}
\tilde P(x)=&a_2 \left(\delta_2-\frac{(2+x+1/\delta_2)\tilde M}{\epsilon} \right)\nonumber\\
&+e^{2\delta_2 x}\left(\frac{1}{\delta_2}-1\right)\frac{\tilde M^2}{\mu}
- \frac{\tilde M^2}{\delta_2 \mu}   e^{\delta_2 (1+x)}.
\end{align}

First,  from (\ref{V2dot}),  we need to  choose the tuning parameter $\delta_2>\frac{1+x}{\mu}$. Then, by choosing
\begin{align}
a_2 > \max & \left\{ 2\rho_1^2e^{2\delta_2}, 2\rho_2^2e^{2\delta_2},\right. \nonumber\\
&\left. -\frac{e^{2\delta_2 x}\left(\frac{1}{\delta_2}-1\right)\frac{\tilde M^2}{\mu}
- \frac{\tilde M^2}{\delta_2 \mu}   e^{\delta_2 (1+x)}}{\delta_2-\frac{(2+x+1/\delta_2)\tilde M}{\epsilon}}  \right\}
\end{align}
to make sure that the matrix $P(x), x\in [0,1]$ is positive definite, we could derive exponential stability of the target error system.
\end{pf}

Then, from the continuity and invertibility of the backstepping transformation (\ref{backtransfo-obs}),
we could derive exponential convergence of the designed observer. Thus, the following theorem is proved.
\begin{thm}
Under the assumptions that  the initial data are in $\left(\mathcal{L}^2([0,1])\right)^3$, the observer (\ref{newvar3-obs}) 
(with the coefficient functions $p_i(x), i=\overline{1,3}$ determined by (\ref{obs-error-kernel})-(\ref{gain-obs}) and with 
the boundary condition $(\ref{bcond-obs})$) exponentially convergent to the system (\ref{newvar3}) in the $\mathcal{L}^2$ sense.
\end{thm}

\section{Output Feedback Control}\label{output}
The controller (\ref{controldef}) requires a full state measurement and the observer is designed to
reconstruct the state over the whole spatial domain based on an output measurement $w(t,0)$. Thus, by combining these two, 
we could design an observer-based output feedback controller.
\begin{thm}
Consider the $(u_1, \ u_2, \ w)^T$-system (\ref{newvar3})-(\ref{bcond})	 together  with the $(\hat u_1, \ \hat u_2, \ \hat  w)^T$-observer 
(\ref{newvar3-obs})-(\ref{bcond-obs}).  For a  given initial condition $(u_1^0, \ u_2^0, \ w^0, \ \hat u_1^0, \ \hat u_2^0, \ \hat  w^0)^T \in \left(\mathcal{L}^2([0,\ 1])\right)^6$
and  the control law
\begin{align}
 U(t)=& -\rho_{1}u_{1}(t,1) - \rho_{2}u_{2}(t,1)\nonumber\\
&+ \int_{0}^{1}\Big [ k_1(1,\xi) \hat u_{1}(x,\xi)  \nonumber \\
& +  k_2(1,\xi)\hat u_{2}(x,\xi) +k_{3}(1,\xi) \hat w(1,\xi) \Big]\,d\xi,\label{outputcontrol}
\end{align}
where $k_1$, $k_2$ and $k_3$ satisfy   (\ref{kernelsyst}) with the boundary condition (\ref{kernelbc}), 
the $(u_1, \ u_2,\ w, \ \hat u_1, \ \hat u_2, \ \hat  w)^T$-system is exponentially stable in the sense of  the $\mathcal{L}^2$-norm.
\end{thm}
\begin{pf}
From the definition of the error variable vector ($\ref{error}$), the combined closed-loop $(u_1, \ u_2,\ w, \ \hat u_1, \ \hat u_2, \ \hat  w)^T$-system of  (\ref{newvar3})-(\ref{bcond}), (\ref{newvar3-obs})-(\ref{bcond-obs}) and (\ref{outputcontrol}) is equivalent with the $(\hat u_1, \ \hat u_2,\ \hat w, \ \tilde u_1, \ \tilde u_2, \ \tilde w)^T$-system of (\ref{newvar3-obs})-(\ref{bcond-obs}), (\ref{newvar3-obs-error})-(\ref{bcond-obs-error}) and (\ref{outputcontrol}). In comparison to  the backstepping transformation (\ref{backtransfo2}) and  (\ref{backtransfo1}), the invertible transformation
\begin{align}
&\hat \psi_i (t,x)= \hat u_{\rm{i}}(t,x)  ~{\rm{ for }}~  i=1,\,2 \label{backtransfo-output2}	\\
 &\hat \chi (t,x)= \hat w(t,x)  - \int_{0}^{x} k_{1}(x,\xi)\hat u_1(t,\xi)\,d\xi \nonumber\\
&~~~~~~~~~~- \int_{0}^{x} k_{2}(x,\xi)\hat u_2(t,\xi)\,d\xi\nonumber\\
&~~~~~~~~~~- \int_{0}^{x} k_{3}(x,\xi)\hat w(t,\xi)\,d\xi\label{backtransfo-output1}
\end{align}
and (\ref{backtransfo-obs}) maps the system (\ref{newvar3-obs})-(\ref{bcond-obs}) into a $(\hat \psi_1, \ \hat \psi_2, \ \hat \chi, \ \tilde \pi_1, $ $\ \tilde \pi_2, \ \tilde \phi)^T$-system, of which the exponential stability can be proved through the following Lyapunov function:
\begin{align}
V(t)=& \int_0^1 a_1 e^{-\delta_1 x}\left(\frac{\hat\psi^2_{1}(t,x)}{\gamma_1} + \frac {\hat\psi^2_{2}(t,x)}{\gamma_2}\right)dx\nonumber\\
&+ \int_0^1  \frac{1+x  }{\mu}\hat\chi^2(t,x) dx\nonumber\\
&+b\left[\int_0^1 a_2 e^{-\delta_2 x}\left(\frac{\tilde  \pi^2_{1}(t,x)}{\gamma_1} + \frac {\tilde \pi^2_{2}(t,x)}{\gamma_2}\right)dx\right.\nonumber\\
&\left.+ \int_0^1\frac{e^{\delta_2 x}}{\mu} \tilde  \phi^2(t,x)dx\right].\label{dotV}
\end{align}
Exponential stability of the $(u_1, \ u_2,\ w, \ \hat u_1, \ \hat u_2,$  $\hat  w)^T$-system is thus proved.
\end{pf}
\section{ Numerical simulations}\label{simulations}
 This section is devoted to the numerical simulations of system (\ref{newvar1}) subject to the boundary conditions
 (\ref{bcond}) using respectively  the  controller $U(t)$  defined in (\ref{controldef}) and (\ref{backtransfo-output1}). 
 Our goal is to demonstrate the performance of the suggested controllers  (\ref{controldef}) and (\ref{backtransfo-output1}) 
 to stabilize system (\ref{newvar1})  around the zero equilibrium.
 For the sake of completeness, we give a short description of the used numerical schemes.
 We employ an accurate finite volume scheme  to advance in time and space
 the hyperbolic evolutionary system (\ref{newvar1}). 
 Elsewhere, for the implementation of the control law  (\ref{backtransfo-output1}),
 a resolution of the kernel PDE's system (\ref{kernelsyst})-(\ref{kernelbc}) on $\mathbb{T}$ is requested.
  For this end, in sight of the triangular shape of   $\mathbb{T}$,  the finite element  setups  are used.
 The solution of the kernel problem are computed accurately by using the quadratic finite  element pair $P_2$.

The mesh of the triangular domain contains   $7655$  degrees of freedom.
 For the evolution equation, the computational  domain  is  the segment $[0,1]$  and is divided uniformly in $100$ cells.
Further,   the CFL number is fixed at $0.9$ and the time step $\Delta t$ in this numerical simulation is given by a  CFL
(Courant-Friedrichs-Lewy) stability conditions
The  initial bottom topography is defined as
${B(0,x)=0.4\Big(1+0.25\exp\Big(-\frac{(x-0.5)^2}{0.003}\Big)\Big),}$
with a gaussian distribution centered at the middle of the domain. 
The initial water level and  its velocity field are computed,
respectively as
$H(0,x)=2.5-B(0,x)$ and $H(0,x)V(0,x)=10sin(\pi x).$\\
From the physical variables of $H(0,x)$, $V(0,x)$ and $B(0,x)$,
the  initial data of the characteristic variables $v$, $u_1$ and $u_2$ are computed  combining
relations (\ref{eq:psi-k}) and (\ref{eq:psi-to-xi}).
 It is interesting to mention that these initial conditions  imply a strong perturbation in the domain.
\subsection{State feedback under subcritical flow regime $(F_r<1)$}
\begin{figure}[!h]
  \centering
\includegraphics[width=0.35\textwidth]{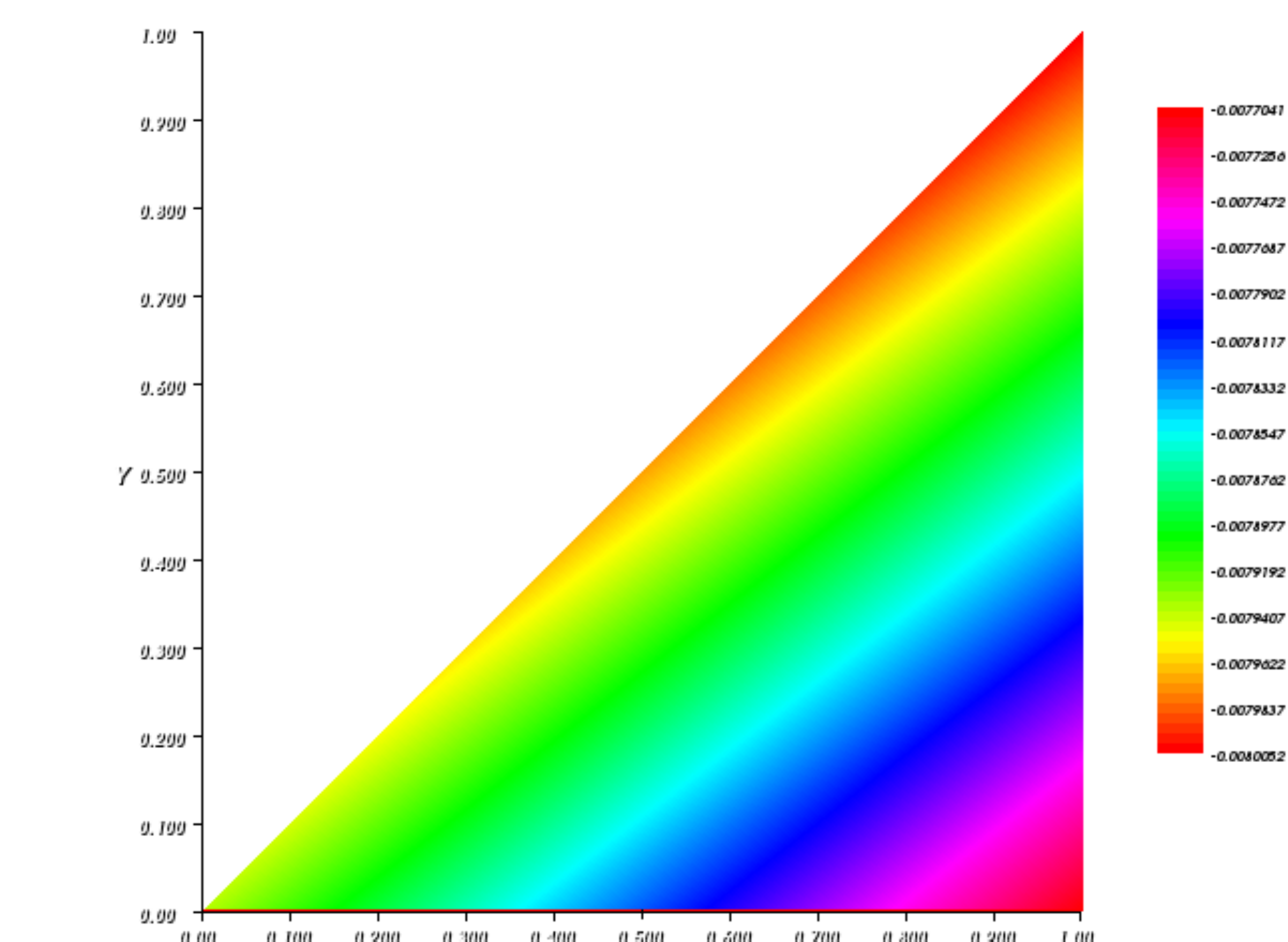}
      \caption{Numerical solution of the kernel  component $k_1$ on  $\mathbb{T}$ 
		 \label{fig:kernel-solution-state}}
\end{figure}

Let us consider the set point ($H^*$, $V^*$, $B^*$) listed in Table \ref{tab:datadetails11-state} (see Appendix)
and make use of the state feedback controller $U(t)$ defined in (\ref{controldef}). The chosen set point leads
to the following characteristic speeds:
$$ \lambda_1 = -1.42,\, \lambda_2= 0.76 \text{ and } \lambda_3 = 7.42.$$
Besides, the Froude number is $Fr=0.6$ which correspond to a subcritical flow regime.
The  coefficients $\alpha_i$, $\theta_i$ and the matrix $\bf \sigma$ are computed 
with the help the characteristics speeds $\lambda_i$. 
Hence  the kernel PDEs (\ref{kernelsyst})-(\ref{kernelbc}) is solved numerically and the value of the kernel 
$k_1$, $k_2$ and $k_3$ at $x=1$  are employed for the  implementation of  the state  
feedback controller (\ref{controldef}). An illustration of the kernel solution $k_1$ is presented in 
Figure \ref{fig:kernel-solution-state}.

\begin{figure}[!h]
  \centering
     \subfigure[Output control law]{\label{fig:controlut-state}\includegraphics[width=0.35\textwidth]{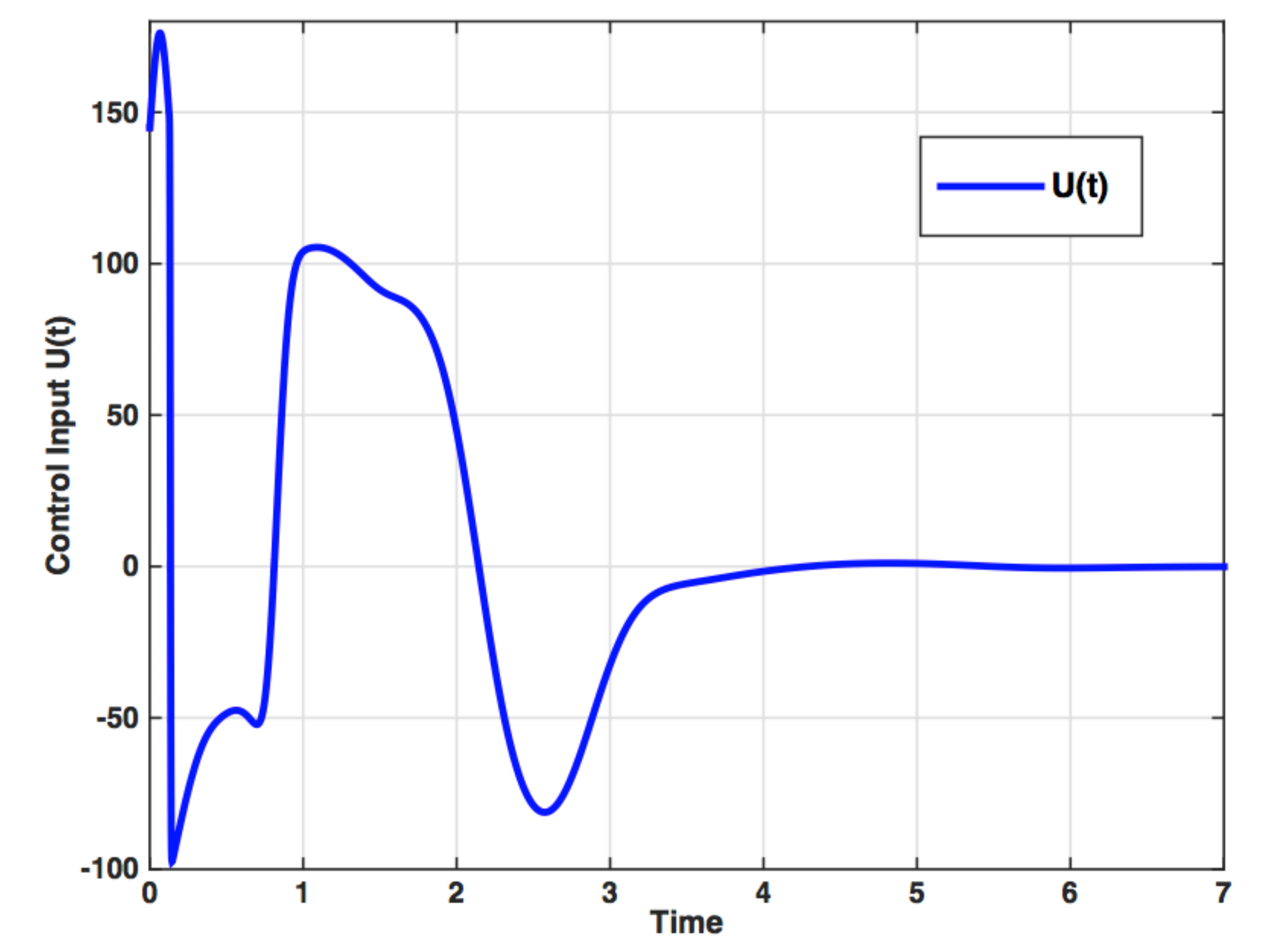}}
     \subfigure[Measured output]{\label{fig:outputyt-state}\includegraphics[width=0.35\textwidth]{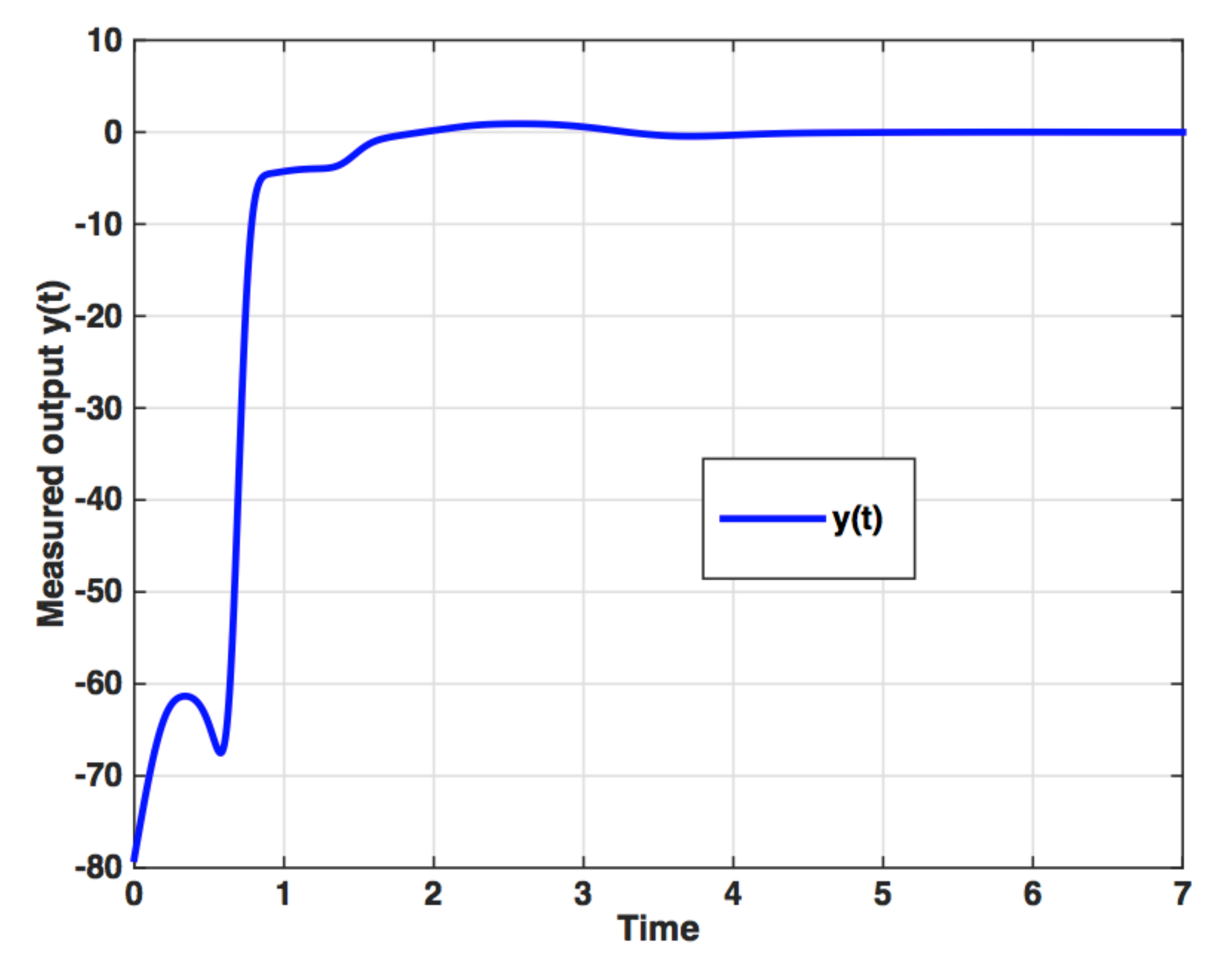}}
      \caption{Evolution in time of the  control law $U(t)$ and the measured output $y(t)$.}
	 \label{fig:control-output-state}
\end{figure}

In Figure \ref{fig:control-output-state} are depicted the behavior in time of the control $U(t)$ and the output measurement $y(t)$. 
Clearly, despite the initial amplitude of $U(t)$, this latter one decreases in time and vanishes after $t\geq4\, s$. 
Let us remind that the implementation of $U(t)$ requires a full-state measurement.
Moreover, likewise $U(t)$ the output measurement $y(t)$ shows the same trend with 
its amplitude decreasing in time and tending to zero after $t\geq 3s$ as can be seen in  Figure \ref{fig:outputyt-state}.

\begin{figure}[!h]
\centering
\includegraphics[width=0.35\textwidth]{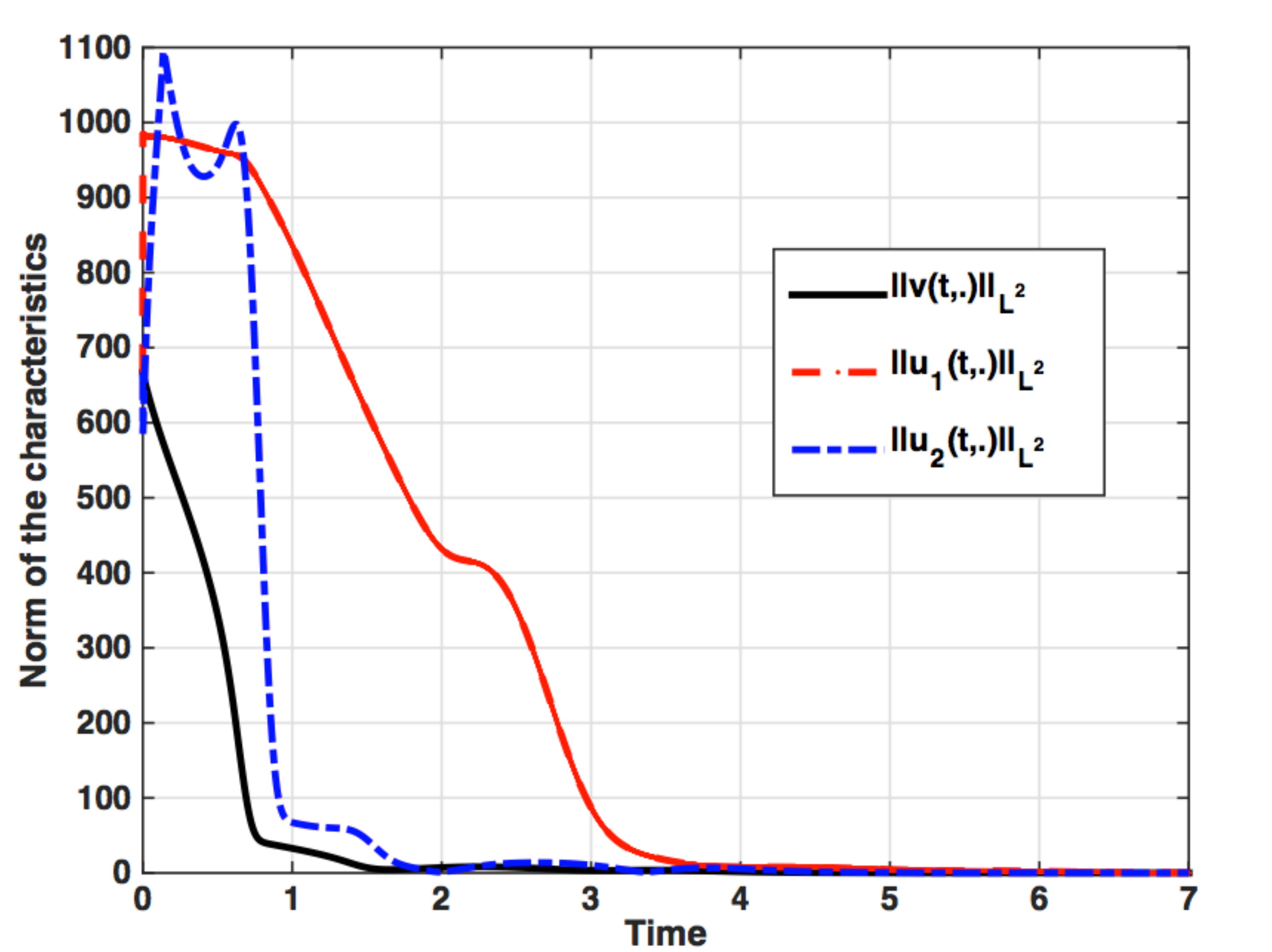}
\caption{Evolution in time of the norm of the characteristic solution.}
\label{fig:normsolution-state}
\end{figure}

Plus, in figure (\ref{fig:normsolution-state}) we plot the evolution  in time of the $\mathcal{L}^2$-norm of the characteristics. As expected from  
the theoritical part we observe that the norm of the characteristics converge to zero. As a result this shows that the system  
(\ref{newvar3}) converge to the zero equilibrium. Thereby the physical linearized model  (\ref{eql_1} ) also converges to  
($H^*$, $V^*$, $B^*$).

\begin{figure}[h]
\centering
\subfigure[Evolution of $u_1(t,x)$]{\label{fig:u1tx-state}\includegraphics[width=0.39\textwidth]{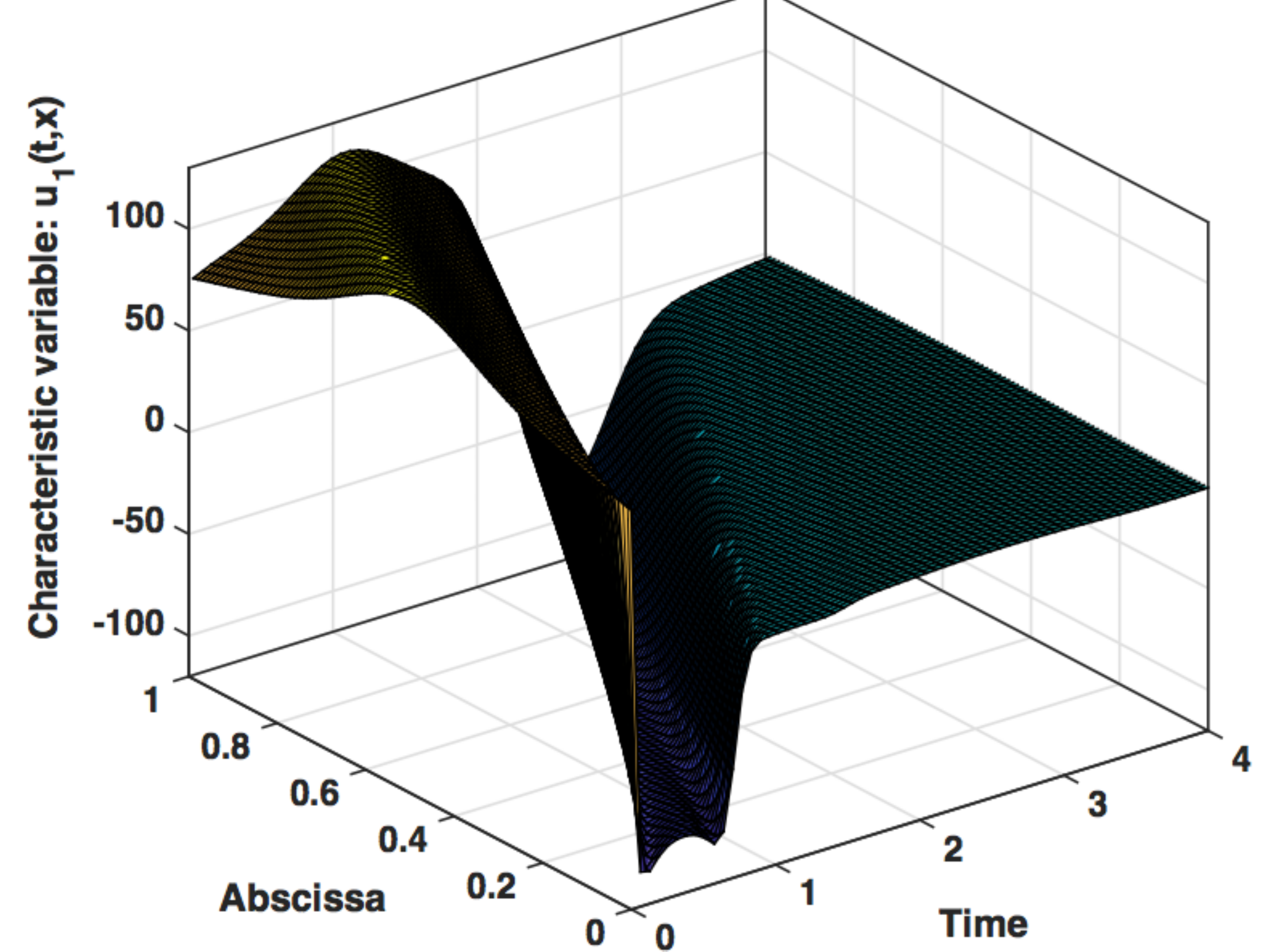}}
\subfigure[Evolution of $u_2(t,x)$ ]{\label{fig:u2tx-state}\includegraphics[width=0.39\textwidth]{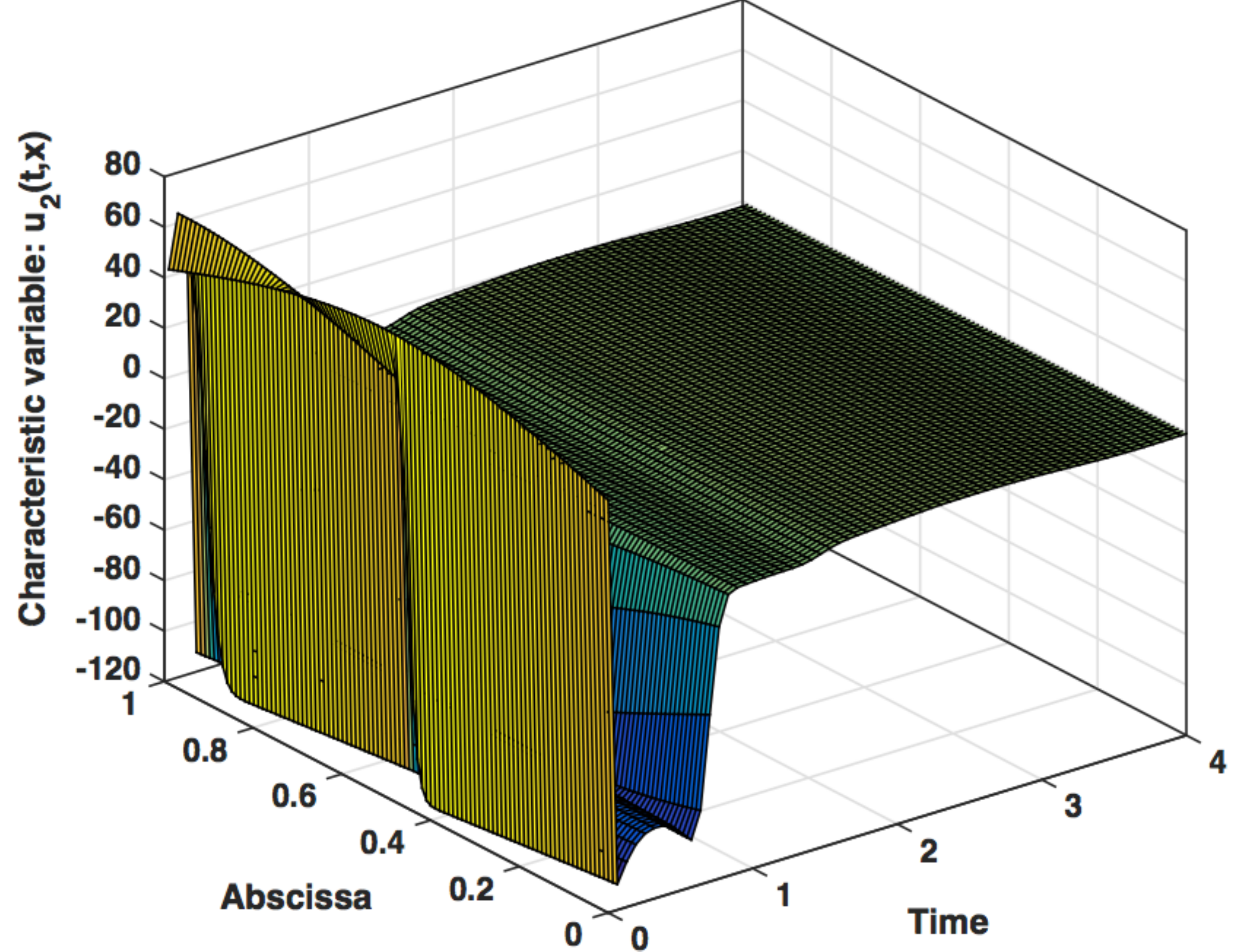}}
\subfigure[Evolution of $v(t,x)$ ]{\label{fig:vtx-state}\includegraphics[width=0.39\textwidth]{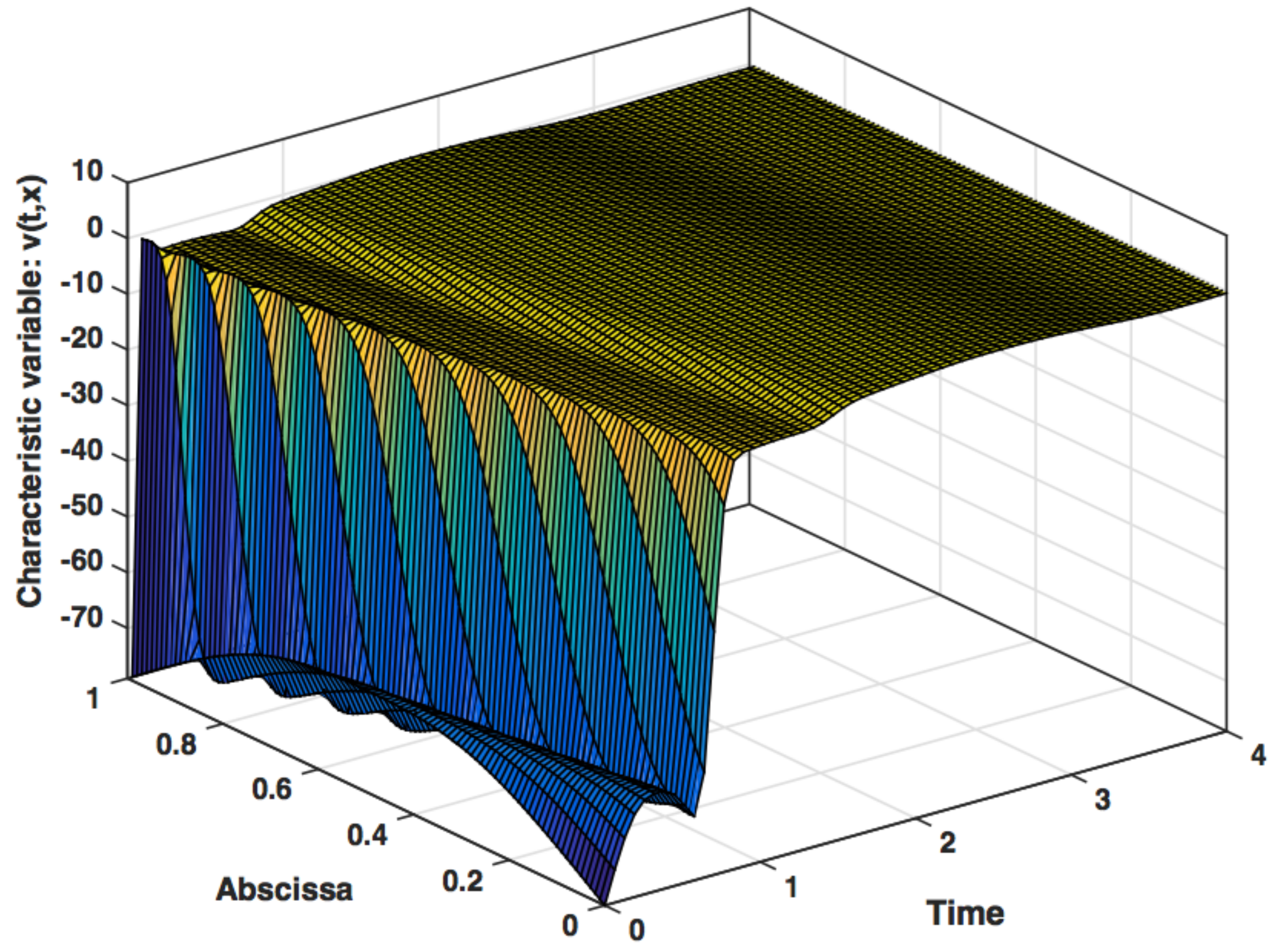}}
\caption{Behavior in time and space of the characteristic solutions.}
\label{fig:abscissa-time-state}
\end{figure}

In addition, Figure \ref{fig:abscissa-time-state} describes the space and time dynamics of the plant
and is consistent with the numerical results presented above.
As time increases, we notice that the perturbation in the overall system decreases and
vanishes later.

\subsection{Output feedback  under supercritical flow regime $(F_r>1)$}
\begin{figure}[htp]
  \centering
   \includegraphics[width=0.35\textwidth]{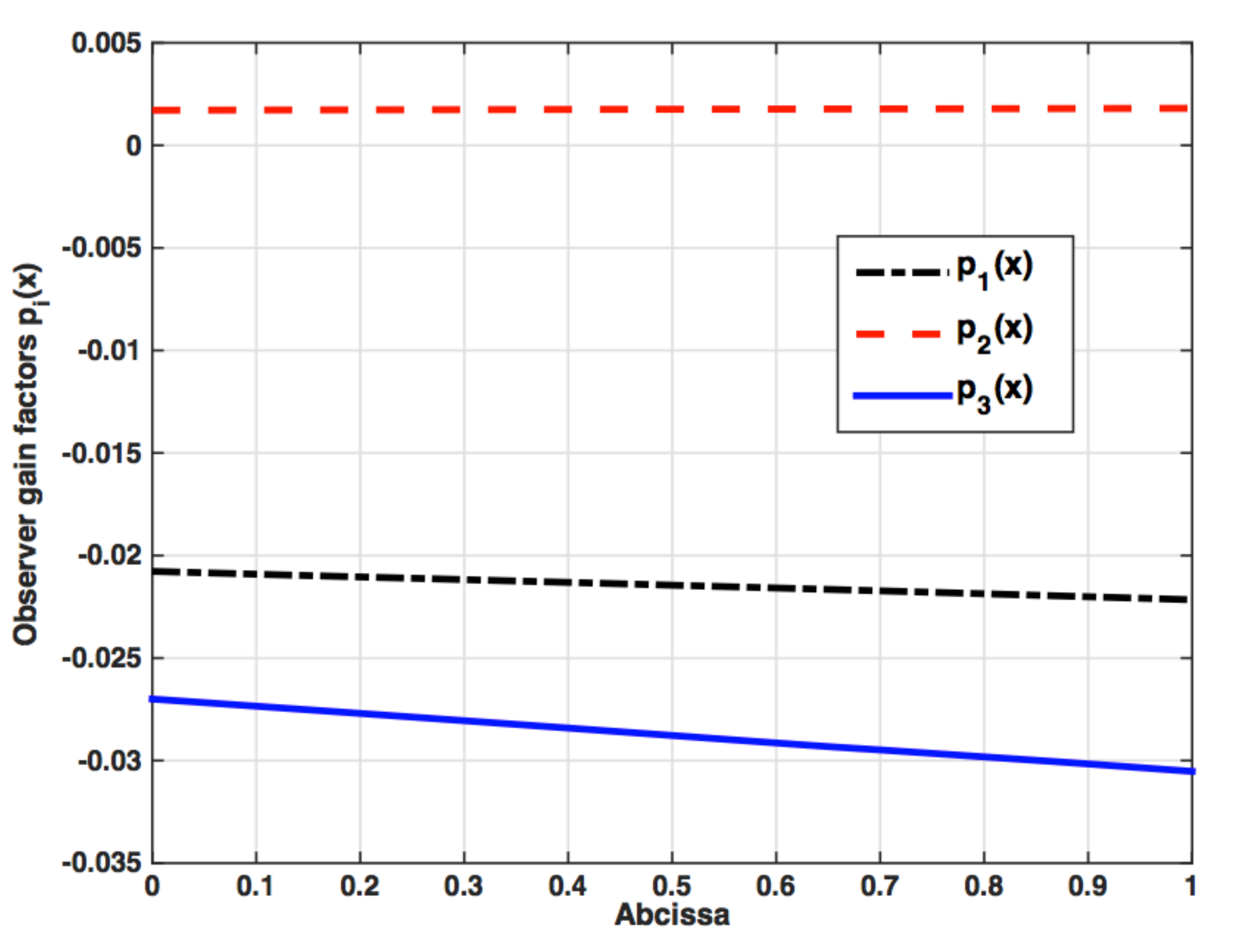}
      \caption{Computed observer gains $p_{\rm i}(x)$.}
	 \label{fig:obsgain}
\end{figure}

 Here, all   parameters of the physical  model are listed in the following Table \ref{tab:datadetails11} given 
 in the Appendix. In this subsection, the dynamic of  the closed-loop system (\ref{newvar1}) together with 
 the output feedback control law (\ref{backtransfo-output1}) is simulated
The set point ($H^*$, $V^*$, $B^*$)   leads to the following  characteristic velocities\\
$$ \lambda_1 = 1.87,\,  \lambda_2=-0.74  \text{ and } \lambda_3 = 8.13.$$
The Froude number is set to $Fr=1.6$. This test case is particularly challenging since the flow 
regime is supercritical. As previously we solve the kernel problem
 As previously the kernel problem (\ref{kernelsyst})-(\ref{kernelbc}) is solved numerically and the
solution is used for the computation of the feedback control law (\ref{backtransfo-output1}).
Not only that, the system (\ref{obs-error-kernel})-(\ref{bc-obs-error-kernel}) is also solved using the  
finite element setup and used to compute the kernel gain $p_{\rm i}(x)$ defined in $(\ref{gain-obs})$. 
This observer gain is represented in Figure \ref{fig:obsgain}.

\begin{figure}[H]
  \centering
  { \includegraphics[height=0.3\textwidth]{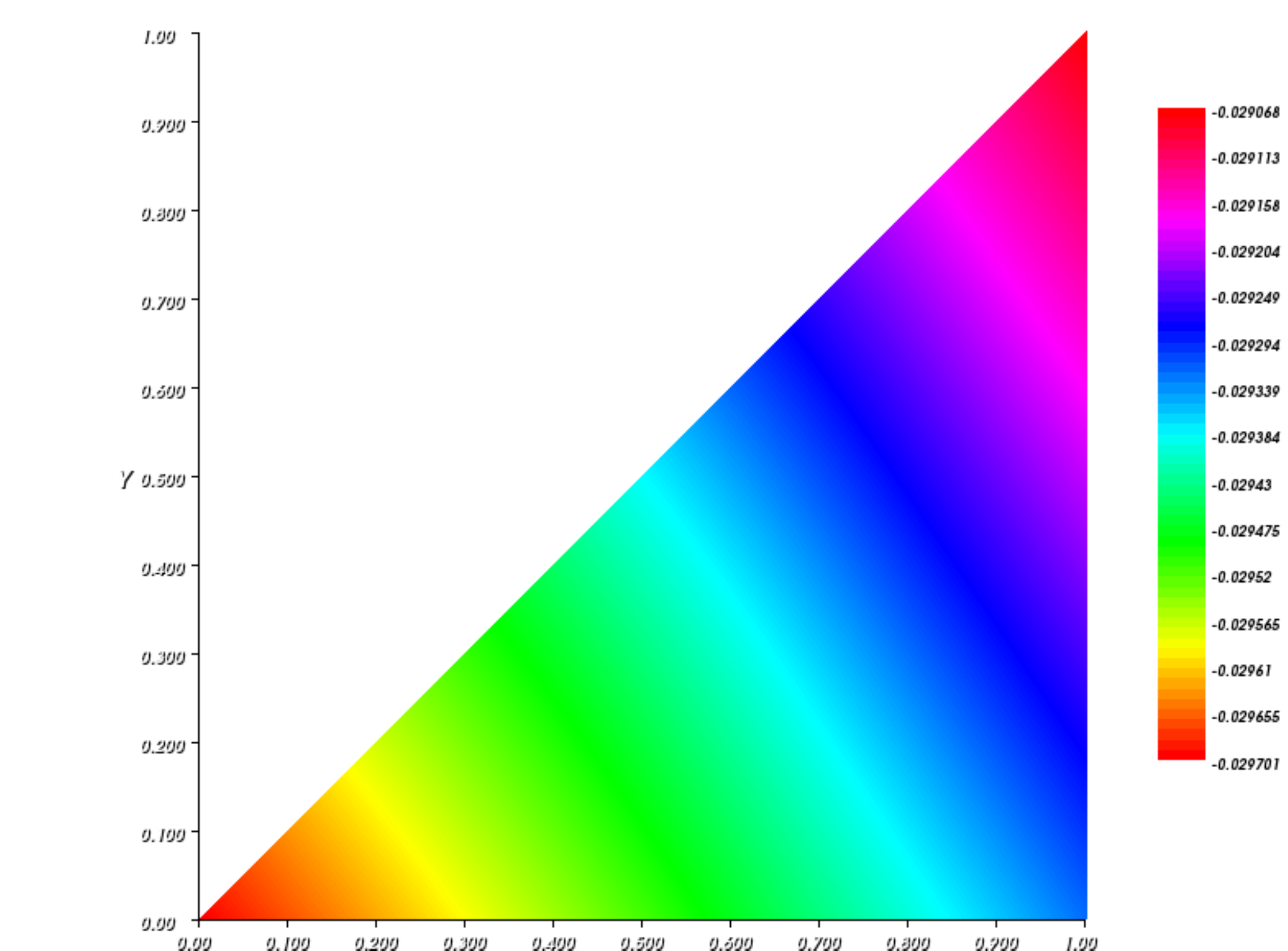}}
      \caption{Solution component $k_1$ on the triangular domain $\mathbb{T}$.}
	 \label{fig:kernel-solution}
\end{figure}

\begin{figure}[H]
  \centering
     \subfigure[Output control law]{\label{fig:controlut}\includegraphics[width=0.35\textwidth]{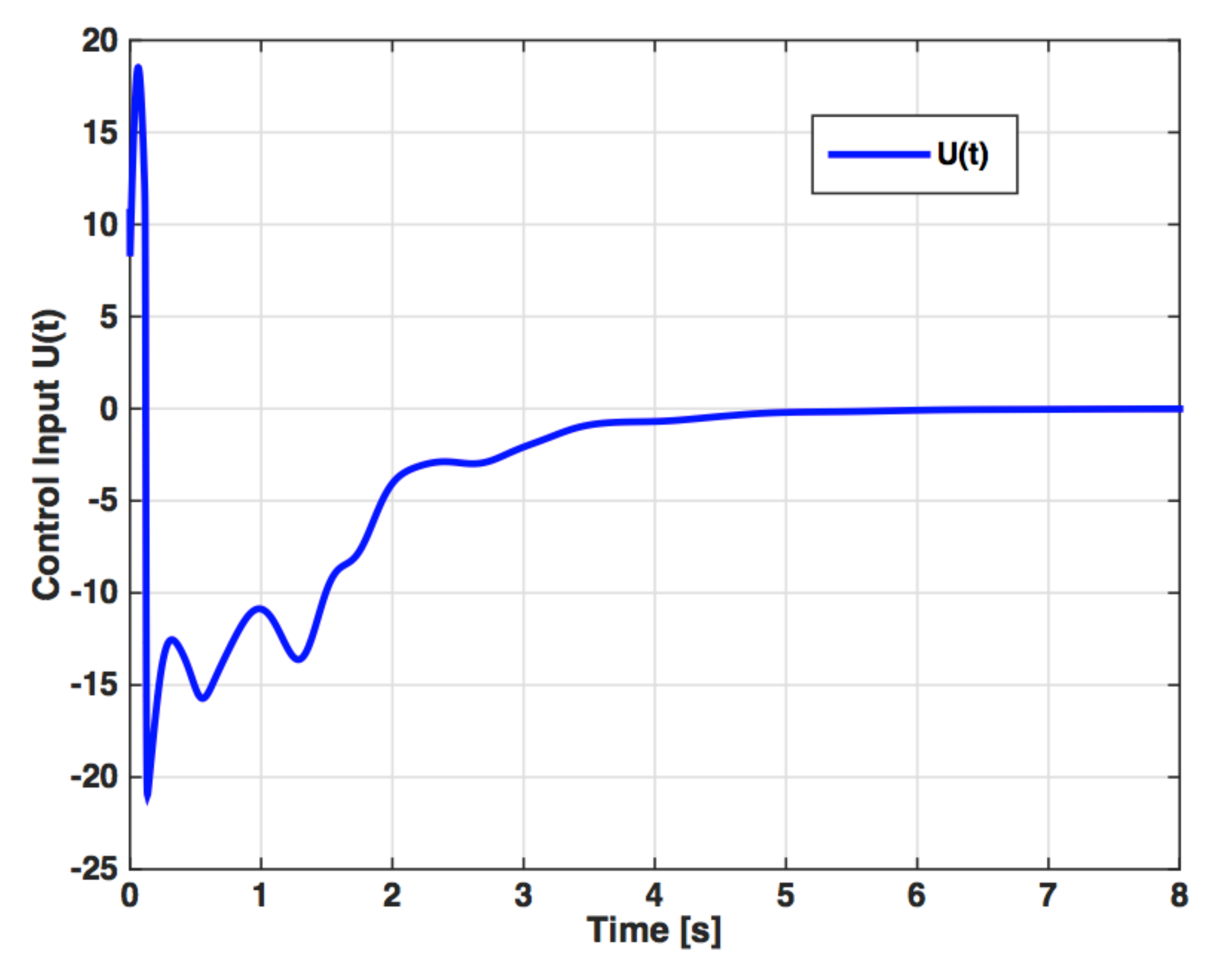}}
      \subfigure[Measured output]{\label{fig:outputyt}\includegraphics[width=0.35\textwidth]{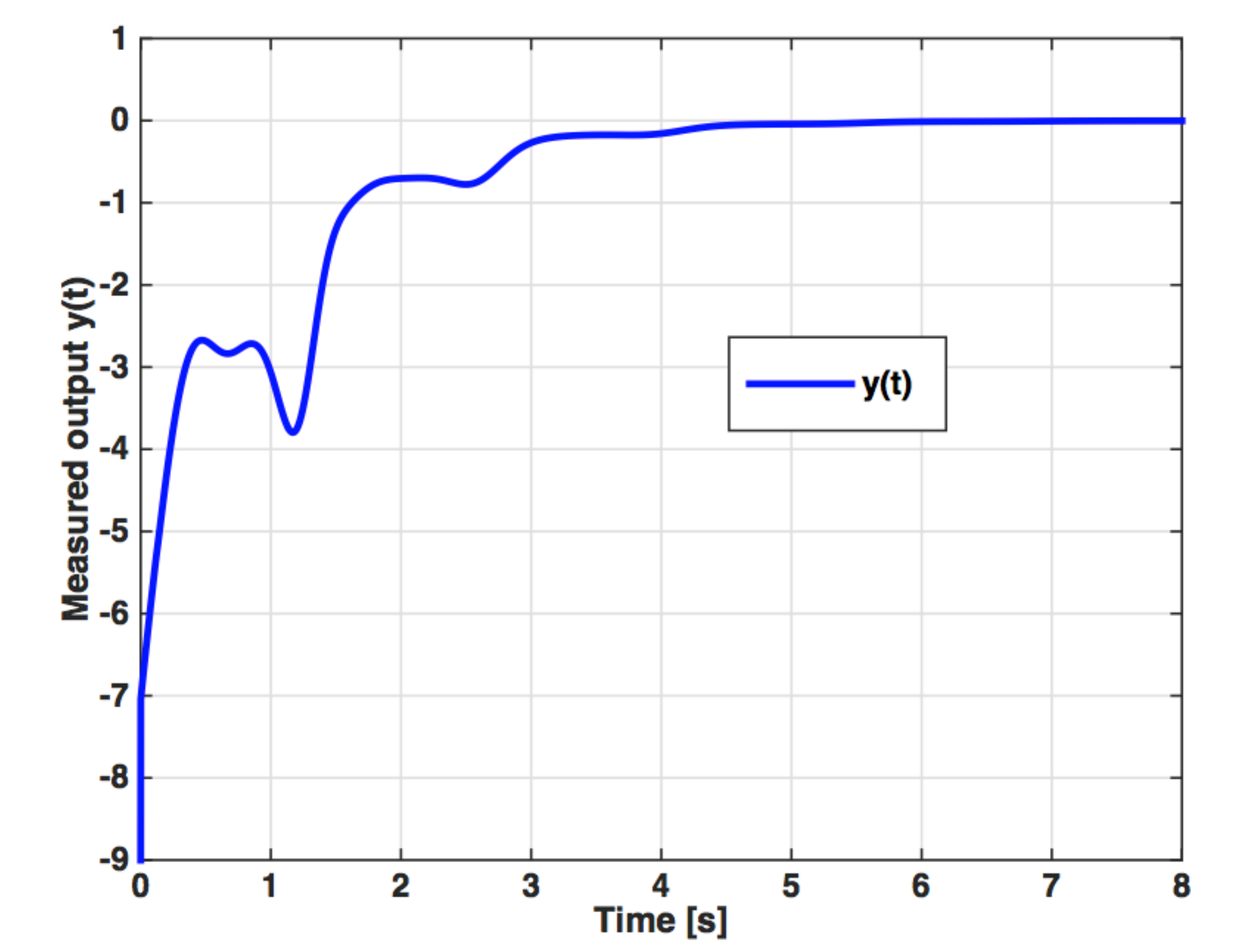}}\hfill
      \caption{Evolution in time of the  control law and the measured output }
	 \label{fig:control-output}
\end{figure}

\begin{figure}[htp]
  \centering
    \subfigure[Ouput feedback control through the backstepping design.]
    {\label{fig:controlback}\includegraphics[width=0.35\textwidth]{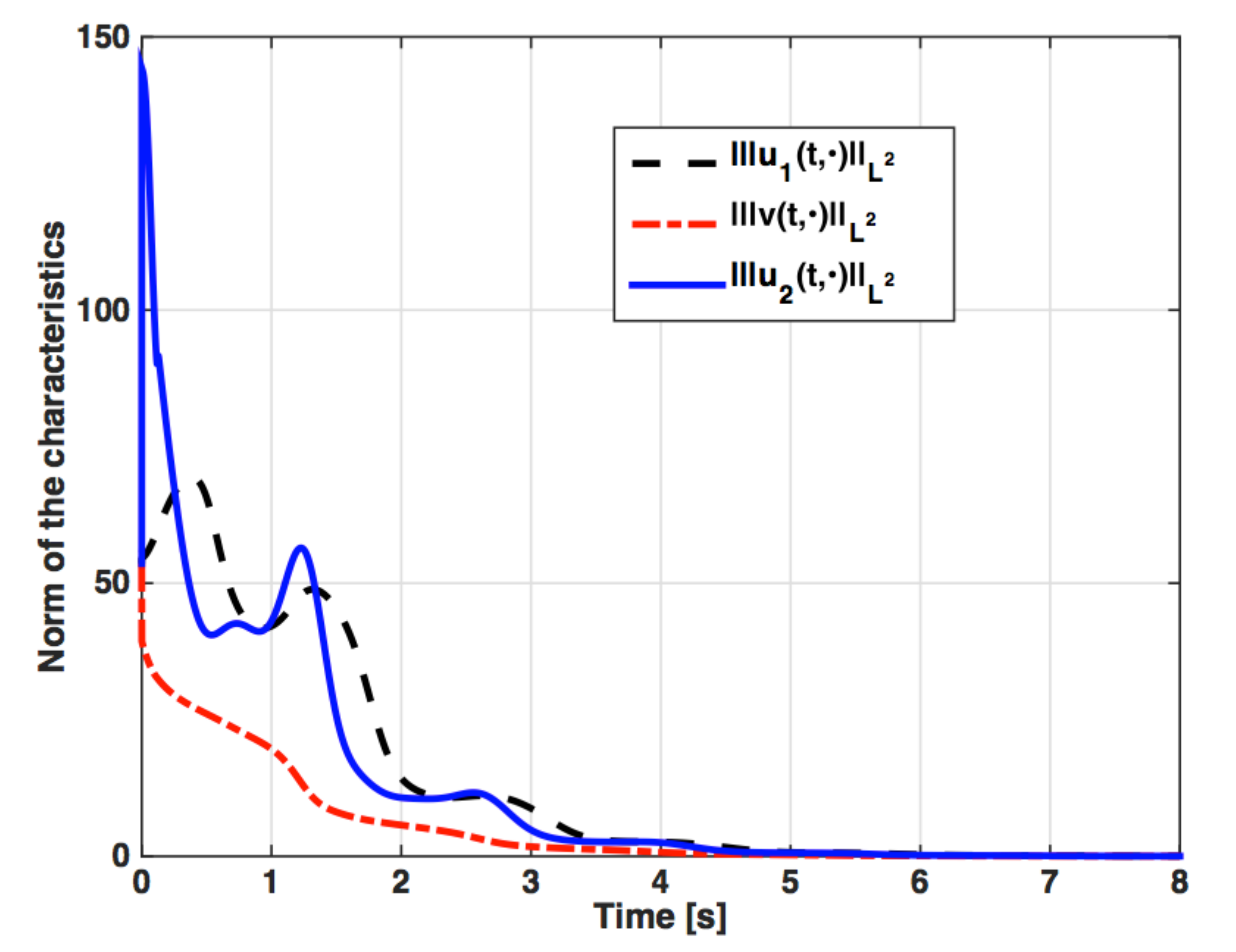}}
     \subfigure[Lyapunov design when the requirements of Theorem $2$ in \cite{Diagne:2011} are not fulfilled.]
     {\label{fig:controllyap}\includegraphics[width=0.35\textwidth]{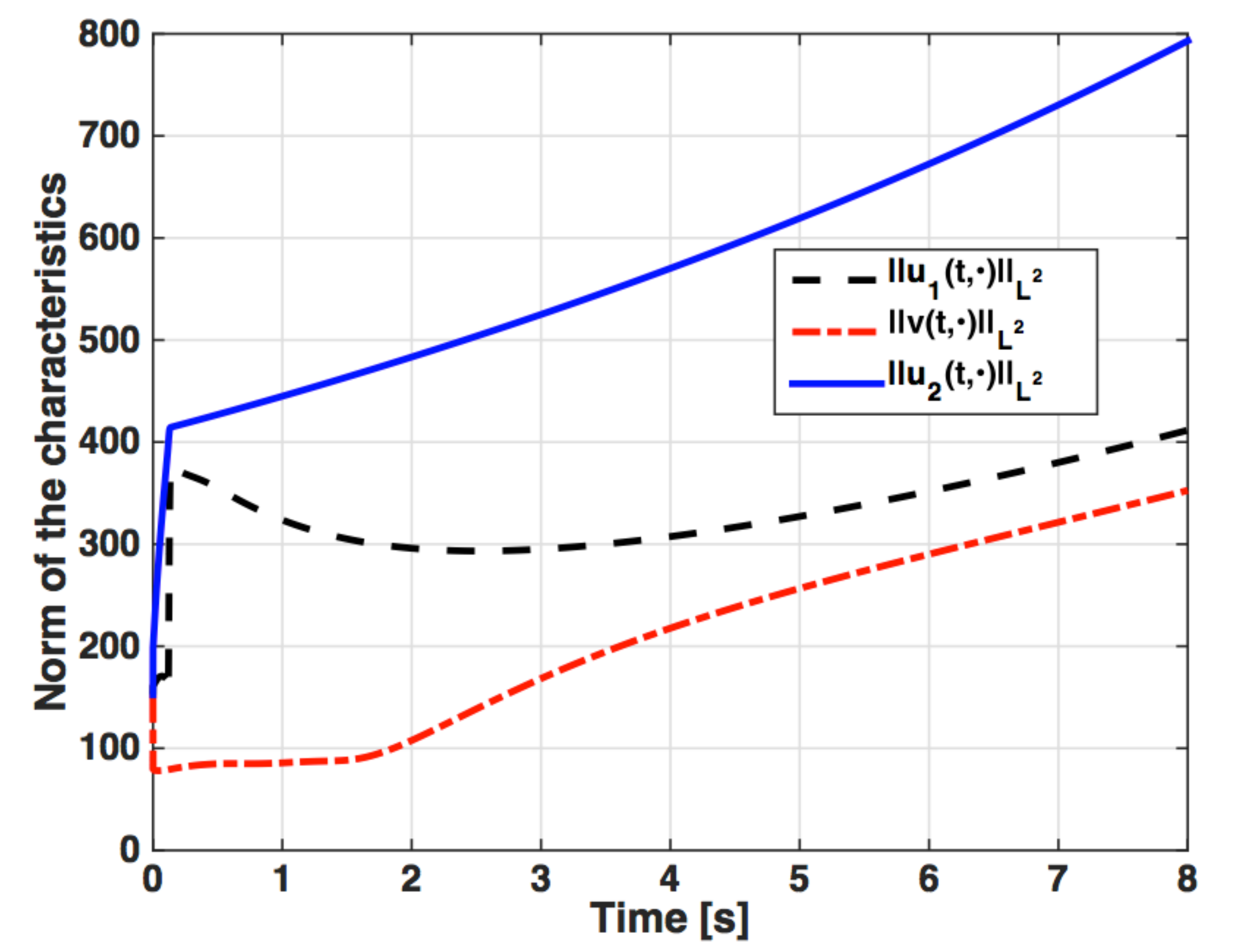}}
      \caption{Evolution in time of the norm of the characteristic solution.}
	 \label{fig:normsolution}
\end{figure}

\begin{figure}[]
  \centering
     \subfigure[Evolution of $u_1(t,x)$]{\label{fig:u1tx}\includegraphics[width=0.35\textwidth]{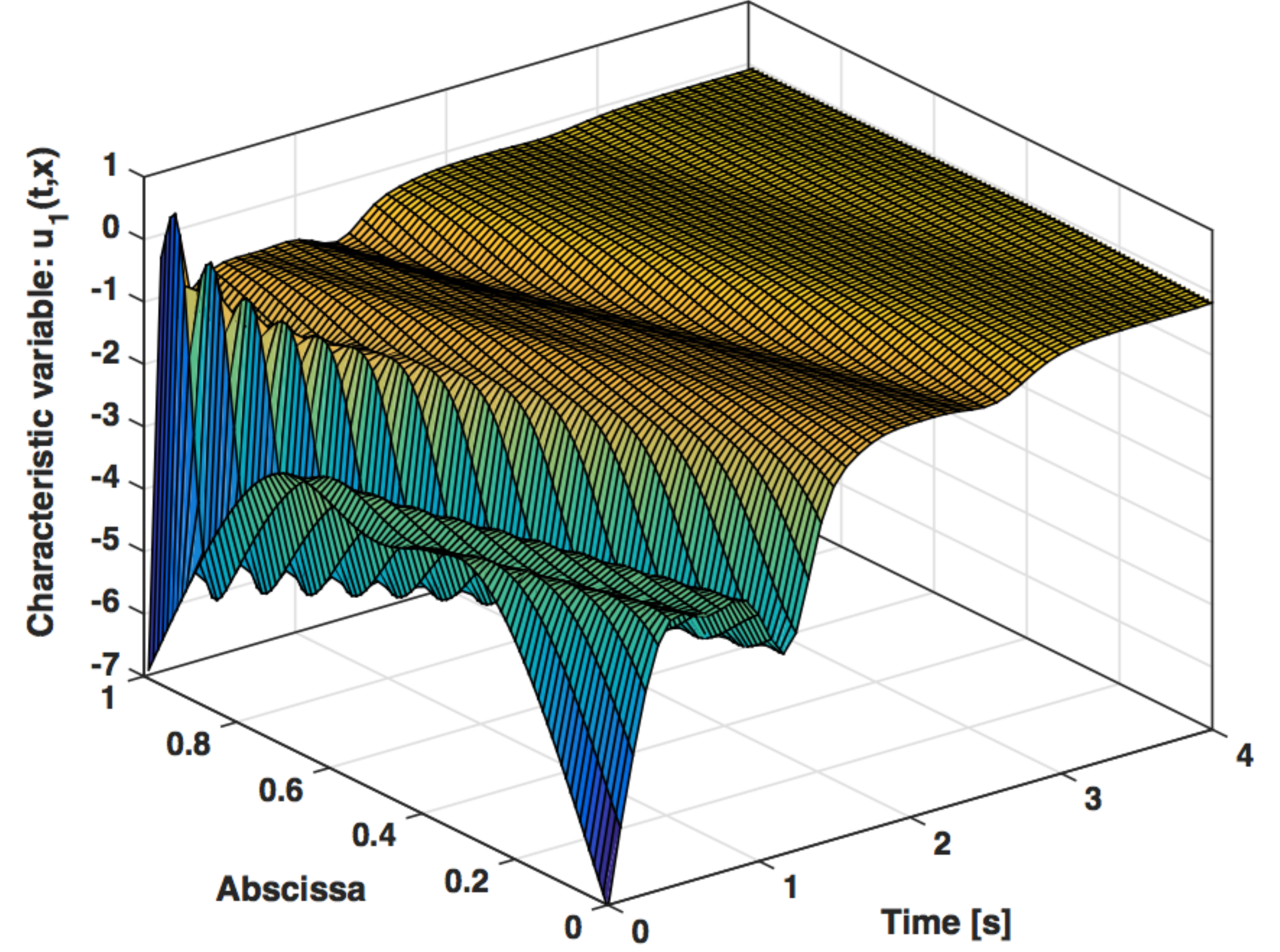}}
     \subfigure[Evolution of $u_2(t,x)$]{\label{fig:u2tx}\includegraphics[width=0.35\textwidth]{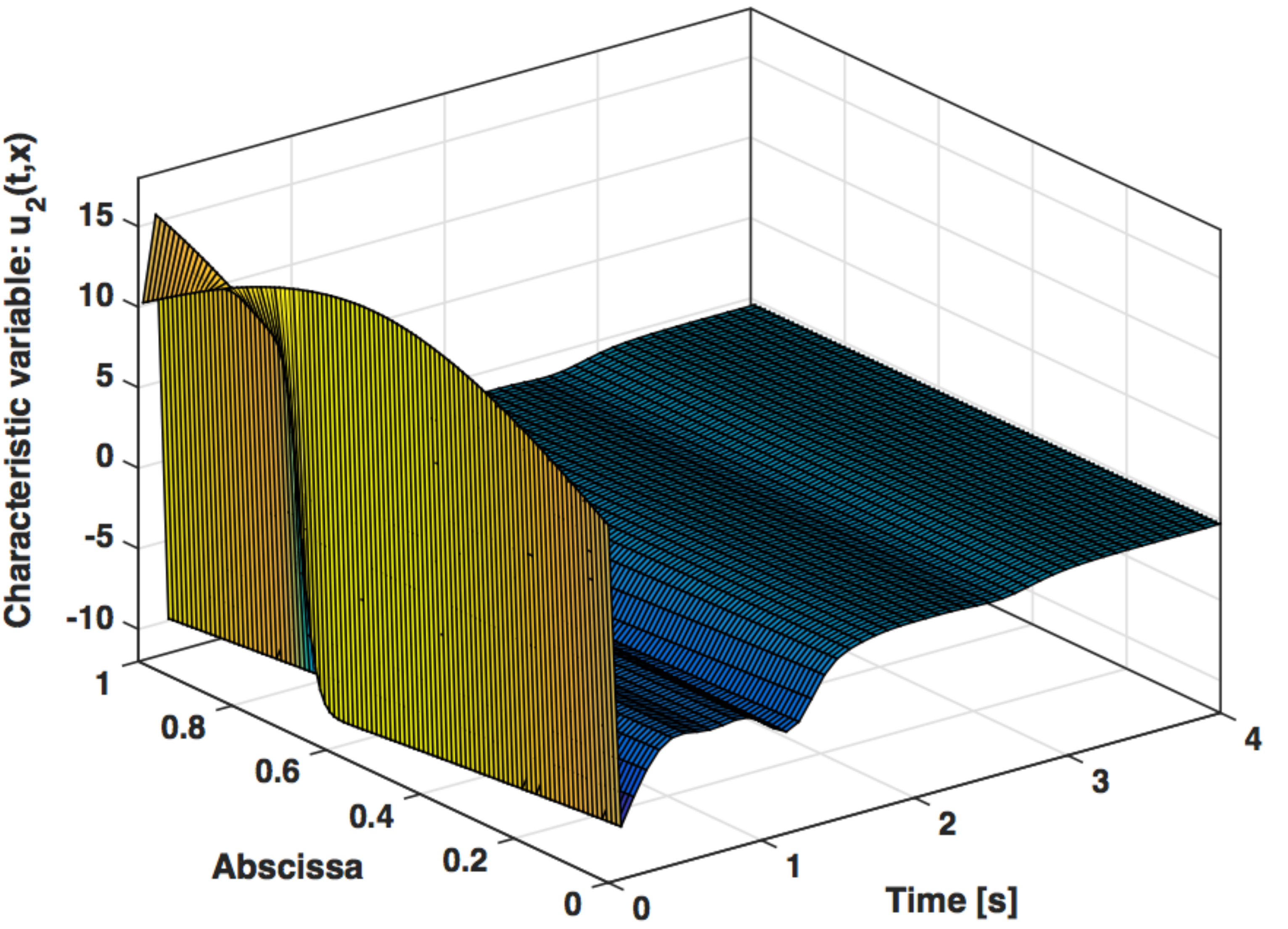}}
     \subfigure[Evolution of $v(t,x)$    ]{\label{fig:vtx}\includegraphics[width=0.35\textwidth]{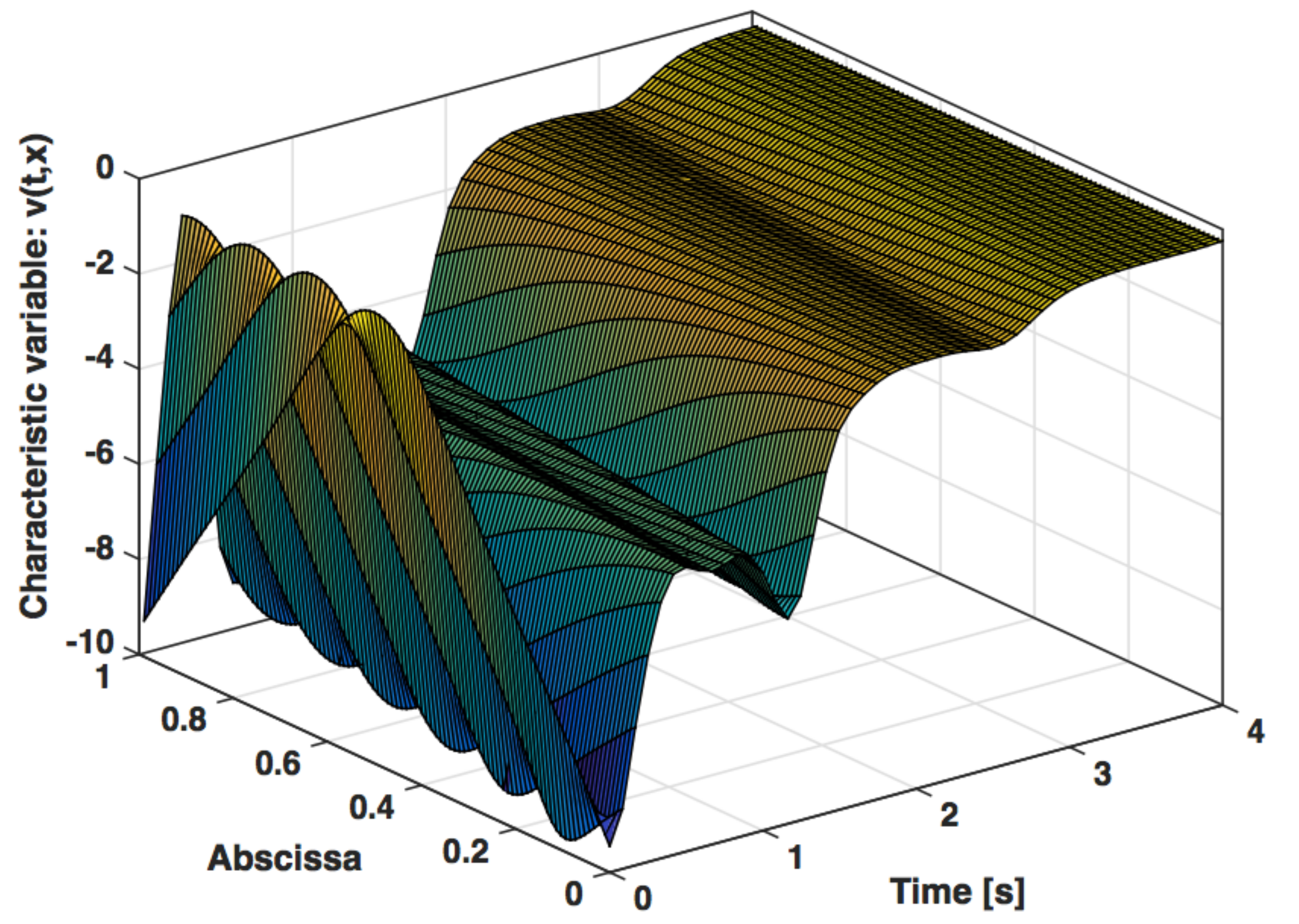}}
      \caption{Behavior in time and space of the distributed states.}
	 \label{fig:abscissa-time}
\end{figure}

In Figure (\ref{fig:kernel-solution}) is depicted  a snapshot of the numerical solution $k_1$  of the kernel PDEs
(\ref{kernelsyst})-(\ref{kernelbc}) with the $x$ and $y$ coordinates  defined as the horizontal and the vertical axis, 
respectively. 
The value of the kernel $k_1$, $k_2$ and $k_3$  at  $x=1$  are  the gain of the designed output
feedback controller (\ref{backtransfo-output1}).
Elsewhere, the computation of the control law (\ref{backtransfo-output1}) requires also the knowledge of the observer. 
Then system (\ref{targetsys-obs})-(\ref{bctarget-obs}) is solved on time and space. 
 Figure  \ref{fig:control-output} shows  the evolution in time of the control input $U(t)$ at downstream  and the output
 measurement $y(t)$ at upstream.  Clearly,  the amplitude  of  $U(t)$  decreases in time and vanishes  for  $t\geq4\, s$.
Moreover, the output measurement $y(t)$ shows  the same trend  with its  amplitude decreasing  in time and tends
to zero after $t\geq 3$.\\
The dynamic of the  $\mathcal{L}^2$-norm is directly related to the magnitude of the propagation speeds
$\lambda_i$ as illustrated in   Figure  \ref{fig:normsolution}. 
Furthermore, we give a comparison of our output feedback law  (Figure \ref{fig:controlback}) to the approach in \cite{Diagne:2011}
(Figure \ref{fig:controllyap}) under this supercritical flow regime (fast rapid flow) where the conditions  of
Theorem 2 (cf \cite{Diagne:2011}) are not fulfilled.
Altogether, our approach exhibits a successful stabilization of system  (\ref{newvar1}) around the zero equilibrium while
instabilities  are noticed  when using the strategy presented in \cite{Diagne:2011}.
Figure \ref{fig:abscissa-time}  describes  the space and  time  dynamics of the plant
and is consistent with the numerical results presented above.
As time increases, we notice that the perturbation in the overall system decreases and
vanishes later.

As can be seen from these numerical simulations, the system (\ref{newvar1})  subject to the feedback control $U(t)$ 
is stabilized around the zero equilibrium as expected from the theoretical part.
We clearly observed the expected qualitative and physical behavior of the designed control regardless the nature of the flow. \\
 This class of state feedback law and output feedback controller we apply here allows to stabilize the SVE system 
 in a optimal time when comparing with the  results presented in \cite{Diagne:2012} where a boundary 
measurement scenario is adopted.


\section{Conclusion and future work}\label{conclusion}

In this paper, a linearized  Saint-Venant Exner model is analyzed for control purposes. The model  describes the evolution of the water  flow coupled with  the  transport of a sediment layer in an open channel. A backstepping state feedback controller, located at the downstream gate of the channel, is first designed for the (exponential) stabilization of the water level and the bathymetry at a desired  equilibrium  set, under the subcritical or supercritical flow regime. Then, based on an exponentially convergent Luenberger observer that  reconstructs  the full state, we design a a backstepping output feedback controller with the measurements at upstream. This controller also achieves the exponential stability of the  linearized  SVE model, for both subcritical and supercritical flow regime.
 Although the backstepping approach offers a  more complicated design than the method developed in \cite{Diagne:2011}, it enables the exponential stabilization of the SVE system without any restriction on the system  and 
 the nature of the flow. It also reduces the number of actuators of  the system: we only need a single boundary control in this paper, but on-line measurements of the water levels at both ends of the spatial domain are needed in \cite{Diagne:2011}. Moreover, simulation results in comparison to \cite{Diagne:2011} are provided to verify that the proposed controller moves beyond those limitations of  \cite{Diagne:2011}.

We emphasize that practically,  such systems are subjected to several types of perturbations and model uncertainties. Thus, an effective control action must take into account of these factors.   We refer the readers to the recent results  proposed in \cite{tang2014a,tang2014b}  on the stabilization of  hyperbolic PDEs with   matched disturbances  at the boundary input. In these papers, the authors employ sliding mode control and active disturbance rejection control to deal with them. Our future objective is to consider robustness issues for this application.  Also, the extension of this approach to a network of flow and sediment transportation  remains   an interesting open problem with a high potential  in real applications.  Among the others,
 the problem of proving local stability of the nonlinear SVE plant under the linear feedback will be very interested to study.
 
 \section*{Appendix}\label{ap}
\begin{itemize}
\item Subcritical flow regime state feedback
\begin{table}[H]
	\centering
\begin{tabular}{|c|c|c|c|c|c|c|c|c|c|c|c|c|c|c|} 
   \hline 
$T$ & $\Delta x$   & $CFL$   & $A $  & $ p$& $C_f$ & $ \rho_1$ & $\rho_2$    \\
    \hline   \hline
     8    &     0.01   &    0.95     &  0.008    & 0.002   &   0.1  &   1.5 &   1.5    \\
    \hline
\end{tabular}
\begin{tabular}{|c|c|c|c|c|c|c|c|c|c|c|c|c|c|c|} 
   \hline 
 $q_1$ & $q_2$ & $H^*$ & $U^*$  & $B^*$   \\
    \hline   \hline
  1   & 1.2   & 2 & 3   & 0.4  \\
    \hline
\end{tabular}
\caption{ Physical parameters and dimensionless numbers}
 \label{tab:datadetails11-state}
\end{table}
\item Supercritical flow regime output  feedback  
\begin{table}[H]
	\centering
\begin{tabular}{|c|c|c|c|c|c|c|c|c|c|c|c|c|c|c|}
   \hline
$T$ & $\Delta x$   & $CFL$   & $A_g $  & $ p_g$& $C_f$ & $ \rho_1$ & $\rho_2$    \\
    \hline   \hline
     8    &     0.01   &    0.9     &  0.003    & 0.002   &   0.1  &   1 &   1.5    \\
    \hline
\end{tabular}
\begin{tabular}{|c|c|c|c|c|c|c|c|c|c|c|c|c|c|c|}
   \hline
 $q_1$ & $q_2$ & $H^*$ & $U^*$  & $B^*$   \\
    \hline   \hline
  1   & 1.2   & 1 & 5   & 0.4  \\
    \hline
\end{tabular}
\caption{ Physical parameters and dimensionless numbers}
 \label{tab:datadetails11}
\end{table}
\end{itemize}

\section*{Acknowledgment}
\addcontentsline{toc}{section}{Acknowledgment}
The first author was supported by grants from \textit{Lisa and Carl-Gustav Esseen} foundation.


\end{document}